\documentclass{article}
\usepackage{amsfonts, amsmath, amssymb, color, a4wide}
\usepackage{hyperref}

\newcommand{\E}{\operatorname{\mathbb{E}}}
\renewcommand{\P}{\operatorname{\mathbb{P}}}
\newcommand{\argmin}{\mathop{\mathrm{argmin}}}
\newcommand{\argmax}{\mathop{\mathrm{argmax}}}
\newcommand{\tr}{\mathrm{Tr}}
\newcommand{\pa}[1]{\left(#1\right)}
\newcommand{\ac}[1]{\left\{#1\right\}}
\newcommand{\cro}[1]{\left[#1\right]}

\newcommand{\eps}{\varepsilon}
\newcommand{\C}{\mathcal{C}}
\newcommand{\cT}{\mathcal{T}}
\newcommand{\cL}{\mathcal{L}}
\newcommand{\ag}{\gtrsim}
\newcommand{\al}{\lesssim}

\newcommand{\X}{\mathbf{X}}
\newcommand{\R}{\mathbb{R}}
\newcommand{\1}{\mathbf{1}}

\newtheorem{thm}{Theorem}

\newtheorem{lem}{Lemma}
\newtheorem{cor}{Corollary}

\begin{document}
\title{Partial recovery bounds for clustering with the relaxed $K$-means}
\author{Christophe Giraud\footnote{Laboratoire de Math\'ematiques d'Orsay, Univ.\ Paris-Sud, CNRS, Universit\'e Paris-Saclay, 91405 Orsay, France. Christophe.Giraud@math.u-psud.fr} and Nicolas Verzelen\footnote{INRA, Montpellier SupAgro, MISTEA,
Univ. Montpellier, France. Nicolas.Verzelen@inra.fr} }

\maketitle

\abstract{
We investigate the clustering performances of the relaxed $K$-means in the setting of sub-Gaussian Mixture Model (sGMM) and Stochastic Block Model (SBM).  After identifying the appropriate signal-to-noise ratio (SNR), we prove that the misclassification error decays exponentially fast with respect to this SNR. These partial recovery bounds for the relaxed $K$-means  improve upon results currently known in the sGMM setting.  In the SBM setting, applying the relaxed $K$-means SDP allows us to handle general connection probabilities whereas other SDPs investigated in the literature are restricted to the (dis-)assortative case (where within group probabilities are larger than between group probabilities). Again, this partial recovery bound complements the state-of-the-art results. All together, these results put forward the versatility of the relaxed $K$-means.}

\section{Introduction}

The problem of clustering is that  of grouping similar ''objects'' in a data set. It encompasses many different instances such as partitioning points in a metric space, or partitioning the nodes of a graph.

\subsection{$K$-means and a convex relaxation}

When these objects can be represented as vectors in a Euclidean space, some of the most standard clustering approaches are based on the minimization of the $K$-means criterion~\cite{Lloyd}. 
Observing $n$ objects and writing $X_{a}\in \R^p$ for the object $a\in \{1,\ldots, n\}$,  the $K$-means criterion of a partition $G=(G_1,\ldots,G_k)$ of $\{1,\ldots, n\}$ is defined as 
\begin{equation}\label{eq:crit_K_means}
 \mathrm{Crit}(G)= \sum_{k=1}^K\sum_{a\in G_k} \bigg\|X_a  - \frac{1}{|G_k|}\sum_{b\in G_k} X_b \bigg\|^2\,, 
\end{equation}
where $\|.\|$ stands for the Euclidean norm. This criterion quantifies the dispersion of each group around its centroid in order to favor homogeneous partitions. 
A $K$-means procedure then aims at finding a partition $\widehat{G}$ that minimizes,  at least locally, the $K$-means criterion. However, solving this problem is  NP-hard and it is even hard to approximate~\cite{awasthi2015hardness}.

In general, iterative procedures such as Llyod's algorithm~\cite{Lloyd} and its variants~\cite{Kmeans++} are only shown to converge to a local minimum of the $K$-means criterion. Alternatively, 
Peng and Wei~\cite{PengWei07}  have suggested to relax the $K$-means criterion to a Semi-Definite Program (SDP) followed by a rounding step. See the next section for a definition.  The resulting program is provably solvable in polynomial time. 
This work is dedicated to promoting Peng and Wei's procedure and some of its variants by (i) putting forward its versatility by handling both vector and general graph clustering problems  and (ii) assessing its near-optimal performances.

\subsection{SubGaussian Mixture Models (sGMM) and Stochastic Block Models (SBM)}

In the computer-science and statistical literature, the most popular approach to assess the performances of a procedure  is the 'model-based' strategy. It assumes there exists a true unknown partition $G$ of the 'objects' and that the data have been randomly generated from a probability distribution rendering this partition. Then, one can assess the performances of  a clustering procedure by comparing the partition estimated from the data  to $G$.

For vector clustering, it is classical to assume that the vectors $X_a$ are distributed according to a SubGaussian Mixture Model (sGMM). In a sGMM with partition $G$, the random variables $X_a$ are assumed to be independent and for $a\in G_k$, the random variable $X_a$ is assumed to follow a subGaussian distribution centered at  $\mu_k\in \mathbb{R}^p$ and with covariance matrix $\Sigma_k$. In other words, variables $X_a$ whose indices $a$ belongs to the same group are identically distributed and  variables $X_a$ and $X_b$ whose indices belong two different groups have different means. See Section \ref{sec:subGauss_definition} for a definition. 

Node clustering in a network has been widely investigated within the framework of Stochastic Block Models (SBM) \cite{holland1983stochastic} and its variants. According to a SBM with partition $G$, the network edges are sampled independently and the probability of presence of an edge between any two nodes $a\in G_k$ and $b\in G_l$ is equal to some quantity $P_{kl}\in [0,1]$ only depending on the groups. In other words, two nodes $a$ and $b$ belonging to the same group in $G$ share the same probability of being connected to any other node $c$.

These two random models have attracted a lot of attention in the last decade. See e.g.~\cite{AbbeReview2017,Moore2017} for two recent reviews on SBM and~\cite{VEMPALA2004,2016arXiv160609190C, mixon2016, LuZhou2016, MartinNIPS,li2017birds,Regev2017,DBLP:journals/corr/abs-1711-07465, DBLP:journals/corr/abs-1711-07211, DBLP:journals/corr/abs-1711-07454} for recent contributions on sGMM. A large body of the literature on these two models focuses on pinpointing the right scaling between the model parameters allowing to recover the partition $G$ from the data. For sGMM, this translates into identifying the minimal distance $\min_{k\neq l}\|\mu_k-\mu_l\|$ within the mixtures means, such that, there exists a clustering procedure, if  possible running in polynomial time, that recovers $G$ with high probability. Most of the works concentrate on two types of recovery: perfect recovery, where one wants to recover exactly the partition $G$ with high probability and weak-recovery where the estimated partition $\widehat{G}$ is only required to be more accurate than random guessing. The goal is then to identify the precise threshold at which perfect or weak recovery can occur. We refer to \cite{AbbeReview2017} for a review of these questions in SBM. Between these two extreme regimes, when the best possible classification is neither perfect nor trivial, the objective is to maximize the proportion of well-classified data. Given two partitions $\widehat G=(\widehat G_{1},\ldots,\widehat G_{K})$ and $G=(G_{1},\ldots,G_{K})$ of $\ac{1,\ldots,n}$ into $K$ non-void groups, we define the proportion of non-matching points  
\begin{equation}\label{misclassify:error}
err(\widehat G,G)= \min_{\pi \in \mathcal{S}_{K}} {1\over 2n} \sum_{k=1}^K \left|G_{k} \bigtriangleup \widehat G_{\pi(k)}\right|,
\end{equation}
where $A\bigtriangleup B$ represents the symmetric difference between the two sets $A$ and $B$ and  $\mathcal{S}_{K}$ represents the set of permutations on $\ac{1,\ldots,K}$. When $\widehat G$ is a partition estimating $G$, we refer to $err(\widehat G,G)$ as the misclassification proportion (or error) of the clustering. The problem of minimizing this error has attracted less attention but see \cite{azizyan2013minimax,guedon2014community, Chin15, YunP14b, AbbeSandon2015a, DAM15, 2016arXiv160609190C, Gao2017, FeiChen2017}  for some related contributions.

Among the polynomial-time clustering procedures, Semi-Definite-Programs (SDP) have proved to be versatile and they have been investigated in a large range of clustering problems, including clustering in SBM \cite{chen2014statistical, guedon2014community, 2015arXiv150705605P, JavanmardE2218, 2016arXiv160206410H, FeiChen2017}, sGMM \cite{2016arXiv160609190C, mixon2016, MartinNIPS,li2017birds} or in block covariance models \cite{pecok,IsingBM}. While not always reaching the exact threshold for weak/perfect clustering in several cases \cite{JavanmardE2218, 2015arXiv150705605P}, SDP algorithms are versatile enough in order to enjoy some robustness properties \cite{2015arXiv150705605P,moitra2016robust,FeiChen2017}, which are not met by more specialized algorithms~(see \cite{moitra2016robust} for more details). However, most SDPs require  the partition to be balanced or that, at least, the size of each group is known in advance. Besides, all SDPs studied for SBM clustering arise as convex relaxations of min-cut optimization  problems \cite{chen2014statistical, guedon2014community, 2015arXiv150705605P, JavanmardE2218, 2016arXiv160206410H, FeiChen2017,li2018convex} and therefore only fall within the framework of assortative SBM where within group probabilities of connection are larger than between group probabilities of connection. In other words, the diagonal entries of $P$ have to be larger than its off-diagonal entries.

\subsection{Our Contribution}

In this work, we provide misclassification error bounds for the relaxed $K$-means of Peng and Wei~\cite{PengWei07} combined with a rounding step, both in the sGMM and the SBM frameworks. 
Compared to other SDPs, this convex relaxation of $K$-means has the nice feature to only require the knowledge of the number of groups (which can sometimes be estimated, as in \cite{pecok}). Hence, there is no need to know the size of the clusters, nor the parameters of the model. The details about this SDP and the subsequent rounding step are given in Section~\ref{sec:def:$K$-means}.

\medskip 

Some of the first partial recovery results for SDPs have been derived using of 
Grothendieck  Inequality~\cite{guedon2014community, 2016arXiv160609190C}. The corresponding  misclassification error bounds scale with a square-root decay with respect to an appropriate signal-to-noise ratio. More recently, Fei and Chen \cite{FeiChen2017} have dramatically improved such bounds in the context of assortative SBM, by proving misclassification error bounds for their SDP that actually exponentially decays with respect to this signal-to-noise ratio. Our results and our proof techniques are inspired by this work.

\medskip  

Let us first give a glimpse of our results on sGMM, by specifying it to the special case of  Gaussian mixture models, with $K$ groups of equal size $m=n/K$ and equal covariance $\Sigma$. The general statement of the results for possibly  unbalanced groups and unequal covariances in sGMM is postponed to  Section~\ref{sec:subGauss}. Write $\Delta=\min_{k \neq l}\|\mu_k-\mu_l\|$ for the minimal Euclidean distance between the means of the components  and write  $R_{\Sigma}=|\Sigma|_{F}^2/|\Sigma|_{op}^2$ for the ratio between the square Frobenius norm of $\Sigma$ and the the square operator norm of $\Sigma$. This ratio can be interpreted as an effective rank of $\Sigma$ and is always smaller than the ambient dimension $p$. In the sequel, $c$ stands for a positive numerical constant. 
Then, Theorem~\ref{thm1} in Section~\ref{sec:subGauss} entails that, with high probability, the proportion of misclassified observations decreases exponentially fast with the signal to noise ratio
\begin{equation}\label{intro:SNR}
s^2={\Delta^2\over |\Sigma|_{op}} \wedge {n\Delta^4\over K |\Sigma|_{F}^2},
\end{equation}
at least, as long as the condition $s^2\geq c K$, or equivalently
\begin{equation}\label{intro:minimal}
\Delta^2 \geq c  |\Sigma|_{op} \pa{1\vee \sqrt{R_{\Sigma}\over n}} K
\end{equation}
is met. The shape (\ref{intro:SNR}) of the signal-to-noise ratio is new and it differs from the classical signal-to-noise ratio $\tilde s^2=\Delta^2/ |\Sigma|_{op}$ considered e.g. in \cite{LuZhou2016}. We explain in Section \ref{sec:iso}  why the exponential decay should be with respect to $s^2$, at least in the isotropic case.
Since $err(\widehat{G},G)\leq 1/n$ implies that the partition $\widehat{G}$ is equal to $G$, the exponential decay with respect to $s^2$ ensures perfect recovery of the clustering with high-probability when $s^2\geq c (K\vee \log(n))$, recovering the results of \cite{MartinNIPS}. It also ensures a better than random guess clustering when (\ref{intro:minimal}) is met, which improves, in high-dimensional setting, upon {state-of-the art} results in \cite{mixon2016,LuZhou2016}. \medskip 

On the SBM side, we explain how the relaxed $K$-means procedure can be applied to general SBM to cluster nodes  presenting similar connectivity profiles. Instead of the previously discussed SDPs that look for a partition  with maximal within-group connectivity, this allows us to handle  general unknown connection matrices $P$ and thereby going  far beyond the assortative case. Denoting by $m$ the size of the smallest group in $G$, we prove that, with high probability, the misclassification proportion decreases exponentially fast with the signal-to-noise ratio
\begin{equation}\label{intro:SNR:SBM}
s^2 = m \cdot\min_{j\neq k}{ \|P_{j:}-P_{k:}\|^2\over |P|_{\infty}},
\end{equation}
at least as long as the condition $s^2\geq cn/m$ is met. Here, $P_{j:}$ stands for the $j$-th row of $P$ and $|P|_{\infty}$ denotes the supremum norm.
Note that this result encompasses sparse graph, where the connection probability may scale as a constant divided by $n$.
When specified to the classical case with  all within-group probabilities equal to $p$ and all between-group probabilities equal to $q$, with $q<p$, and all groups of the same size, we recover the results obtained by \cite{FeiChen2017} for a relaxed version of the MLE,  but  without knowing that we are in the assortative case, nor knowing the group sizes.

\subsection{Connection to the literature}
Only a few papers have previously proved theoretical properties on the relaxed $K$-means of \cite{PengWei07}.  \cite{awasthi2015relax,iguchi2015tightness,li2017birds}
obtain perfect recovery results for the so-called stochastic ball models and Gaussian mixture models, and  \cite{mixon2016} provides bounds on the estimation of the centers of the means in the sGMM, under a condition stronger than (\ref{intro:minimal}). Closer to the present paper, \cite{pecok,MartinNIPS} (see also \cite{cord}) prove perfect recovery results in the setting of block covariance models and sGMM.
To the best of our knowledge, the main result of \cite{MartinNIPS} provides the weaker condition in high-dimension ($p\geq n$) ensuring perfect recovery with polynomial-time algorithm in the sGMM. This condition is $s^2\geq c (K\vee \log(n))$, with $s^2$ defined by (\ref{intro:SNR}).  
 Theorem~\ref{thm1} below extends this result to the partial recovery regime,  in the sense that the main result in \cite{MartinNIPS} can be recovered from Theorem~\ref{thm1}.

\medskip  

In the sGMM setting, the paper \cite{LuZhou2016} derives partial recovery results for Lloyd algorithm (with a suitable initialization). To the best of our knowledge, these results are the strongest ones in the literature. In the setting discussed above, they prove a decay of the misclassification proportion exponentially fast relative to $\tilde s^2=\Delta^2/|\Sigma|_{op}$, with a minimal signal-to-noise requirement 
$$ \Delta^2 \geq c |\Sigma|_{op} K^2 \pa{1\vee {pK\over n}}.$$
 Here, the requirement is proportional to $K^2|\Sigma|_{op}$ in low-dimension, which is larger than our $K|\Sigma|_{op}$ by a factor $K$.
 In high-dimension, since we always have $R_{\Sigma}\leq p$, the factor ${pK/ n}$ is also much larger than our $\sqrt{R_{\Sigma}/n}$ in condition (\ref{intro:minimal}). This more limited range of validity is partially due to the fact that \cite{LuZhou2016} investigates exponential decay with respect to $\tilde s^2=\Delta^2/|\Sigma|_{op}$, rather than  the suitable signal-to-noise ratio $s^2$ given by (\ref{intro:SNR}). We refer to Section~\ref{sec:iso} (and Appendix \ref{sec:bayes}) for a short explanation of this point. Yet, compared to us, \cite{LuZhou2016} have a tight constant in the exponential rate.

During the wrap-up of this paper, we became aware of an independent and simultaneous work of Fei and Chen \cite{FeiChen2018}, which also investigate partial recovery in sGMM with another SDP. They show interesting connections between their SDP and the error of the supervised classification problem with known centers, and derive some partial recovery bounds based on it. Their results also have  a more limited range of validity than ours, as they require groups of the same size and a minimal signal-to-noise condition of the form 
$$ \Delta^2 \geq c |\Sigma|_{op} \pa{K \pa{1\vee {p\over n}}+ \sqrt{Kp\log(n)\over n}}\ ,$$ 
instead of (\ref{intro:minimal}). As before, this more limited condition is partially due to the fact that they investigate exponential decay with respect to $\tilde s^2=\Delta^2/|\Sigma|_{op}$, instead of $s^2$ given by (\ref{intro:SNR}). 

\medskip 

In the SBM setting, most results on partial recovery \cite{Chin15,YunP14b,AbbeSandon2015a,Gao2017,FeiChen2017} cover the assortative setting. The papers \cite{Chin15,YunP14b,Gao2017} investigate some two-steps procedures based on a spectral algorithm. The papers  \cite{YunP14b,Gao2017} derive tight misclassification bounds for their algorithm, showing sharp exponential decay with respect to the signal-to-noise ratio.  Closer to us, \cite{FeiChen2017} proves similar results for an SDP, with less tight constants than \cite{YunP14b,Gao2017}, but a wider range of validity. 
Compared to these results, our results does not provide sharp constants as in \cite{YunP14b,Gao2017}. Yet, they provide some new results for partial recovery in non-assortative cases and they only require a mild condition on the size of the smallest cluster. To the best of our knowledge, (i) our results are the first results about clustering with an SDP in non-assortative cases and (ii)  the only known exponential bounds for partial recovery in general SBM are those of \cite{AbbeSandon2015a} which handle the sparse setting where the matrix $P$ scales as $P=P_{0}/n$, with $P_{0}$ fixed and $n\to \infty$. Their results are optimal in the vicinity of the weak recovery threshold. Our results cover a setting with slightly more signal (the misclassification error has to be  smaller than $e^{-cK}$), and the results do not overlap. In particular, as discussed in Section~\ref{sec:SBM}, our exponential rate (\ref{intro:SNR:SBM}) involved in Theorem~\ref{thm2} is faster by (at least) a factor $K$ than the exponential rate involved in 
\cite{AbbeSandon2015a}, though the rate of \cite{AbbeSandon2015a} cannot be improved in the vicinity of the weak recovery threshold. Hence both results are more complementary than comparable. 
\medskip

Since our work has been made available, two follow-up papers have extended and complemented our results.~\cite{ChenYang2018} have proved an exponential clustering error with respect to the SNR $s^2$ \eqref{intro:SNR} for sGMM in infinite dimensional Hilbert spaces with common covariance matrix. For spherical mixtures with  $K=2$,~\cite{Ndaoud2018} has proved a minimax lower bound showing that the clustering error $\exp(-c s^2)$ is  optimal.

\subsection{Organization and notation}
The paper is organized as follows. In Section~\ref{sec:def:$K$-means}, we recall the relaxed $K$-means SDP derived by \cite{PengWei07} and we explain how the partition is derived from the solution of the SDP. Section~\ref{sec:subGauss} covers the sGMM and Section~\ref{sec:SBM} covers the SBM. We explain the main lines of the proofs and discuss the main arguments in Section~\ref{sec:outline}. Finally, the full proof of the two main theorems can be found in Section~\ref{sec:proof:thm1} and Section~\ref{sec:proof:thm2}. 
\smallskip

\noindent {\bf Notation.} To any matrix $M$ we denote by ${\rm Row}(M)$ the set of its rows, by $|M|_{1}$ the $\ell^1$ norm of its entries, by $|M|_{op}$ its operator norm with respect to the $\ell^2$ norm, by $|M|_{F}$ its Frobenius norm, by $|M|_{*}$ its nuclear norm and by $\tr\pa{M}$ its trace. 
We also associate to a diagonal matrix $D$, the pseudo-norm $|D|_{V}=\max_{a}D_{aa}-\min_{a}D_{aa}$.
Besides, for any  $A,B$ with the same dimensions, we write $\langle A,B\rangle=\sum_{ab} A_{ab}B_{ab}$ for its canonical inner product. For a vector $x$ we write $\|x\|_2$ for its Euclidean norm and $\langle.,.\rangle$ for the corresponding inner product.  In the sequel, $\1$ stands for the indicator function.

For two sequences $u_n$ and $v_n$ (possibly depending on other parameters), we write $u_n \lesssim v_n$ (resp. $u_n \gtrsim v_n$) when there exists some numerical constant $c>0$ such that $u_n \leq c v_n$ (resp. $u_n \geq c v_n$). Given $x,y\in \mathbb{R}$, $x\vee y$ (resp. $x\wedge y$) stands for the maximum (resp. the minimum) of $x$ and $y$.

\section{Relaxed $K$-means}\label{sec:def:$K$-means}
We have $n$ ``objects" that we want to cluster. In the case of sGMM, these objects are $p$-dimensional vectors, and in the case of SBM they corresponds to the the nodes of a graph. For each object $a$, we have a $p$-dimensional vector of observations:  In the sGMM setting, the observation related to $a$ is the vector $X_{a}\in \mathbb{R}^p$ and  in the SBM the observation is the vector  $X_{a}\in \ac{0,1}^n$  recording presence/absence of edges between $a$ and the other nodes (hence $p=n$ is this case).  We denote by $\X\in \R^{n\times p}$, the $n\times p$ matrix  whose $a$-th row is given by $X_{a}$. In particular, in the SBM case, $\X$ is simply the adjacency matrix of the graph. 

Peng and Wei \cite{PengWei07} have observed  that any partition $G$ of $\{1,\ldots, n\}$ can be uniquely represented by a $n\times n$ matrix $B\in \mathbb{R}^{n\times n}$ such that $B_{ab}=0$ if and only if $a$ and $b$ are in different group and  $B_{ab}=1/|G_k|$ if $a$ and $b$ are in the same group $G_k$. The collection of such matrices when $G$ spans the collection of all partitions with $K$ groups may be described as 
$$\mathcal{P}=\ac{B\in \R^{n \times n}:\textrm{symmetric},\ B^2=B,\ \tr(B)=K,\ B1=1,\ B\geq 0}\ . $$
Here, $B\geq 0$ means that all entries of $B$ are nonnegative and $B^2$ refers to the matrix product of $B$ with itself.
Peng and Wei \cite{PengWei07} have shown that minimizing the classical $K$-means criterion~\eqref{eq:crit_K_means} is equivalent to maximizing $\langle \X\X^T, B\rangle$ over the space $\mathcal{P}$. Writing $\tilde{B}$ for such a maximizer, the  $K$-means clustering is obtained from $\tilde B$ by grouping together indices $a,b$ which have a non-zero entry $\tilde B_{ab}$. 

The constraint set $\mathcal{P}$ is non-convex and solving the $K$-means problem is NP-hard  \cite{awasthi2015hardness}. 
Peng and Wei \cite{PengWei07} then propose to relax the constraint set $\mathcal{P}$ by dropping the condition $B^2=B$  to consider
$$\C=\ac{B\in \R^{n \times n}:\ \textrm{Positive Semi Definite,}\ \tr(B)=K,\ B1=1,\ B\geq 0},$$
and hence solve the relaxed $K$-means SDP
\begin{equation}\label{relaxed-$K$-means}
\widehat B \in \argmax_{B\in \mathcal{C}} \langle \X\X^T, B\rangle.
\end{equation}
Obviously, the solution $\widehat{B}$ does not necessarily belong to $\mathcal{P}$ and then does not provide a clustering.  One has therefore to rely on a rounding step to obtain a proper partition. If $\widehat{B}$ is close to the true matrix $B$, one should expect that rows of $\widehat B$ belonging to the same group are similar. This is why the final step is obtained by applying a clustering algorithm on the rows  of $\widehat B$.
 As in \cite{FeiChen2017}, we apply here an approximate $K$-medoids on the rows of $\widehat B$. Let us detail this final step.

To any partition  $G=(G_{1},\ldots,G_{K})$ of $\ac{1,\ldots,n}$ into $K$ (non empty) groups, we can associate a partition matrix $A\in \R^{n\times K}$ defined by $A_{ik}=\1_{i\in G_{k}}$. Let us denote by $\mathcal{A}_{K}$ the set of such matrices and by $\textrm{Rows}(\widehat B)$ the set of the rows of $\widehat B$. Then a $\rho$-approximate $K$-medoids on the rows  of $\widehat B$ is a pair $(\widehat A, \widehat M)$ with $\widehat A \in \mathcal{A}_{K}$, $\widehat M\in\R^{k\times n}$, ${\rm Rows}(\widehat M)\subset {\rm Rows}(\widehat B)$ and fulfilling
\begin{equation}\label{eq:$K$-medoid}
 |\widehat A\widehat M-\widehat B|_{1}\leq \rho \min_{A\in\mathcal{A}_{K},\ {\rm Rows}(M)\subset {\rm Rows}(\widehat B)} |AM-\widehat B|_{1}.
 \end{equation}
We refer to \cite{kmedian} for a polynomial-time algorithm producing such an output $(\widehat A, \widehat M)$ with $\rho=7$. Then a partition is obtained from $\widehat A$ by setting $\widehat G_{k}=\ac{i: \widehat A_{ik}=1}$.

In the sequel, $\widehat{G}$ is said to be a relaxed $K$-means solution if $\widehat{G}$ is derived from any $7$-approximate $K$-medoids on the rows of $\widehat{B}$ obtained in~\eqref{relaxed-$K$-means}.  All our partial recovery bounds are for this partition $\widehat G$.

As a final remark,  $K$-means is known to suffer from a bias which tends to produce groups of similar width, see e.g. \cite{MartinNIPS}. 
As explained in \cite{pecok, MartinNIPS}, $K$-means and its relaxed version can be debiased when useful (e.g. for high-dimensional mixtures with unequal traces $\tr(\Sigma_k)$).
We refer to Section~\ref{sec:hetero} for details.

\section{Clustering sub-Gaussian mixtures}\label{sec:subGauss}
\subsection{Model}\label{sec:subGauss_definition}
We observe $n$ independent random vectors $X_{1},\ldots,X_{n}\in \R^p$. We assume that there exists an unknown partition   $(G_{1},\ldots,G_{K})$ of $\ac{1,\ldots,n}$ and $K$ unknown $p$-dimensional vectors $\mu_{1},\ldots,\mu_{K}\in\R^p$, such that 
$$X_{a}=\mu_{k}+E_{a}\quad\textrm{for any}\ a\in G_{k},$$
with $E_{a}$ centered, independent with covariance $\Sigma_{k}$. We recall that $m=\min_k|G_k|$ stands for the size of the smallest group. 

The larger the Euclidean distance between two centers $\Delta_{jk}=\|\mu_{k}-\mu_{j}\|$, the more easily we can recover the unknown partition from the observations $X_{1},\ldots,X_{n}$. Hence, we denote by $\Delta=\min_{j\neq k} \Delta_{jk}$ the minimal distance between two distinct centers, which will represent the signal part in the signal-to-noise ratio. 

The hardness of the clustering problem also depends on the concentration of the random vectors $E_{a}$ around zero. A common distributional assumption when analyzing clustering is the subgaussiannity assumption. For a centered random vector $Z\in \mathbb{R}^p$ and $L>0$, we say that $Z$ is SubG($LI_p$) if the random variables $Z_1,\ldots,Z_{p}$ are independent and  $\E[\exp(tZ_i)]\leq \exp(t^2 L^2/2)$ for all $t\in \mathbb{R}$ and $i=1,\ldots, p$. For sub-Gaussian mixture, we make the following distributional assumption.
\medskip

\noindent{\bf Assumption A1.} There exists $L>0$ such that $\Sigma_{k}^{-1/2}E_{a}$ is SubG$(L^2I_{p})$.
\medskip

Under this assumption, two quantities mainly drive the noise width in the signal-to-noise ratio: the maximum scaled operator norm of the covariances $\sigma^2=L^2\max_{k}|\Sigma_{k}|_{op}$, and the maximum scaled Frobenius norm of the covariances $\nu^2=L^2\max_{k}|\Sigma_{k}|_{F}$. 

Actually, as shown in Theorem~\ref{thm1} below, the misclassification error of the relaxed $K$-means decreases exponentially fast with the signal-to-noise ratio 
\begin{equation}\label{SNR:subG}
s^2={\Delta^2\over \sigma^2}\wedge  {m\Delta^4\over \nu^4} ,
\end{equation} 
where $m$ denotes the size of the smallest cluster. This particular  definition of the signal-to-noise ratio is new and is further discussed below.
\medskip

\noindent\underline{Remarks:}
\begin{enumerate}
\item When the random variable $E_{a}$ is normally distributed with covariance $\Sigma_k$, it fulfills Assumption A1 with $L=1$.  
\item We observe that the random variable $E_{1},\ldots,E_{n}$ are all sub-Gaussian SubG$(\sigma^2I_{p})$, and we can always upper-bound $\nu^2$ by $\nu^2\leq p \sigma^2$, with equality in the spherical case where the covariances  $\Sigma_{k}$ are proportional to the identity matrix. Yet, this upper bound is crude when the covariances are far from being proportional to the identity matrix.
\end{enumerate}

\subsection{Partial recovery bound}
Let us denote by $\Gamma\in\R^{n\times n}$ the diagonal matrix with entries $\Gamma_{aa}=\tr(\Sigma_{k})$ for $a\in G_{k}$. 
As shown in \cite{MartinNIPS} (see also \cite{pecok}), when the trace of the covariances $\Sigma_{1},\ldots,\Sigma_{K}$ are unequal, i.e. when $\Gamma$ is not proportional to the identity, it is useful to de-bias the relaxed $K$-means (\ref{relaxed-$K$-means}) by removing from $\X\X^T$ a preliminary estimator $\widehat \Gamma$ of $\Gamma$. 
This estimator can be $\widehat \Gamma=0$ {(no correction)} when the covariances have equal traces or be equal to  (\ref{estimator:Gamma}) as  defined in  Section~\ref{sec:hetero} when the trace of the covariances are unequal. This leads to the so-called Pecok estimator \cite{pecok, MartinNIPS}
\begin{equation}\label{pecok}
\widehat B \in \argmax_{B\in \mathcal{C}} \langle \X\X^T-\widehat \Gamma, B\rangle\ .
\end{equation}
As explained in Section \ref{sec:def:$K$-means}, $\widehat{G}$ is then computed as any $7$-approximate $K$-medoid solution of $\widehat{B}$.
We recall the notation $|D|_{V}=\max_{a}D_{aa}-\min_{a}D_{aa}$.

\begin{thm}\label{thm1}
There exist three positive constants $c,c',c''$ such that the following holds. 
Assume that Assumption A1 holds,
\begin{equation}\label{cond:Gamma}
\Delta^2 \geq {64|\Gamma-\widehat\Gamma|_{V}\over m},
\end{equation}
and that $s^2$ (as defined in \eqref{SNR:subG}) satisfies $s^2 \geq c'' n/m$, then
 the proportion of misclassified points is upper bounded by $$err(\widehat G,G)\leq e^{-c's^2}\ ,$$
with probability at least $1-c/n^2$,  
\end{thm}

Leaving aside Condition \eqref{cond:Gamma} that will be discussed below, one observes that the  misclassification error decreases exponentially fast with respect to the signal-to-noise ratio $s^2$ as soon as  $s^2$ is large enough. 
We further discuss this result in the next two {paragraphs}, first discussing the cases of equal covariance traces and then turning to unequal trace case.

\subsection{Equal traces case}\label{sec:iso}
We assume here that all the covariance matrices $\Sigma_{k}$ have equal trace. Hence, $\Gamma$ is proportional to the identity and $|\Gamma|_{V}=0$. So, when choosing 
 $\widehat \Gamma=0$ in \eqref{pecok}, the condition (\ref{cond:Gamma}) becomes $\Delta^2\geq 0$ which always holds.
\smallskip

\noindent {\bf Non-trivial recovery.} Theorem \ref{thm1} ensures non-trivial recovery as soon as $s^2 \ag n/m$. 
Introduce $R_{\Sigma}$ as the ratio 
\begin{equation}\label{eq:definition_R_Sigma}
R_{\Sigma}={\nu^4\over \sigma^4}={\max_{k=1,\ldots,K}|\Sigma_{k}|_{F}^2\over \max_{k=1,\ldots,K}|\Sigma_{k}|^2_{op}}\leq \max_{k=1,\ldots,K} {|\Sigma_{k}|^2_{F}\over |\Sigma_{k}|^2_{op}}\leq p\,,
\end{equation}
which can be interpreted as an effective rank of the mixture model.
In order to compare this result with those in the literature, let us discuss this condition in the special case of balanced partition with $K$ groups of equal size $m=n/K$.
Theorem~\ref{thm1} guaranties  a non-trivial recovery as soon as $s^2\ag K$, or equivalently
\begin{equation}\label{eq:non-trivial}
{\Delta^2\over \sigma^2} \ag \pa{1\vee \sqrt{R_{\Sigma}\over n}} K,
\end{equation}
with a misclassification error upper-bounded by $e^{-c'K}$ with high-probability. Taking for granted that, as we advocate below,  the exponential decay $e^{-c's^2}$ is optimal in some cases, we cannot hope for a weaker condition than (\ref{eq:non-trivial})  to ensure  a misclassification error of at most $e^{-c'K}$ when $K\al \log(n)$ (for larger $K$ a misclassification proportion of $e^{-c'K}\leq 1/n$ ensures perfect recovery). Yet, it is possible to get a weaker dependence in $K$ when $K\ag \log(n)$. Actually, in a large sample setting where $n\ag p^3K^2 \log(pK)$, \cite{VEMPALA2004} derives clustering guaranties for an iterate spectral clustering under the weaker condition $\Delta^2 \ag \sigma^2\big(\sqrt{K\log(n)}+\log(n)\big)$. Nevertheless, we emphasize that our result, contrary to theirs, holds in the high-dimensional regime $p\geq n^{1/3}$. For the different goal of learning the means in a large sample size setting, some recent papers \cite{DBLP:journals/corr/abs-1711-07465, DBLP:journals/corr/abs-1711-07211, DBLP:journals/corr/abs-1711-07454} have shown that a separation $\Delta \geq K^{\epsilon}$ in the number $K$ of cluster is enough for learning the means in polynomial time, when the sample size is larger than  $n\geq poly(p^{1/\epsilon},k)$. For the same question, \cite{mixon2016} has shown that, in a large sample size setting, the relaxed $K$-means succeeds to learn the means when $\Delta^2/\sigma^2 \ag {K}^2$, which is a stronger requirement than (\ref{eq:non-trivial}). Turning back to our problem of deriving exponential bound for the misclassification proportion, \cite{LuZhou2016} provides such bounds for the Lloyd algorithm under the minimal requirement
\begin{equation}\label{LZ:condition}
 {\Delta^2 \over \sigma^2} \gtrsim   K^2 \pa{1\vee {pK\over n}},
\end{equation}
which, again, is stronger than (\ref{eq:non-trivial}).
 To the best of our knowledge, our result is the first result of this kind for an SDP in this setting. We mention yet, that in an independent and simultaneous work, Fei and Chen \cite{FeiChen2018} have derived a similar in spirit result in the very precise setting where the groups are of equal size. Actually, for an SDP taking as input that all groups have the same size $n/K$, \cite{FeiChen2018} shows non-trivial recovery when 
\begin{equation}\label{cond:FC2018}
{\Delta^2\over \sigma^2} \ag \pa{1\vee {p\over n}} K +\sqrt{Kp\log(n)\over n}.
\end{equation}
Since we always have $R_{\Sigma}\leq p$, the requirement (\ref{cond:FC2018}) of \cite{FeiChen2018} or (\ref{LZ:condition}) of \cite{LuZhou2016} are stronger than (\ref{eq:non-trivial}), especially in the practical case where $p$ is larger than $n$, but the effective rank $R_{\Sigma}$ is small compared to $p$.

\medskip

\noindent{\bf Intermediate regime.} In the intermediate regime where $n/m \al s^2 \al (n/m)\vee \log(n)$, the misclassification rate of our procedure decays at the exponential rate $s^2= \tfrac{\Delta^2}{\sigma^2}\wedge \frac{m\Delta^4}{\nu^4} $. To simplify the discussion, we assume again that the clustering is made of $K$ groups of equal size $m=n/K$. For spherical mixtures, the misclassification rate of the Bayes classifier decays   at the exponential rate $\tilde{s}^2= \tfrac{\Delta^2}{\sigma^2}\geq s^2 $. In view of the definition of $s^2$ and $R_{\Sigma}$, Theorem \ref{thm2} ensures that relaxed $K$-means achieves this optimal rate as soon as 
\begin{equation}\label{eq:condition_bayes_rate}
 \frac{\Delta^2}{\sigma^2}\gtrsim \pa{1\vee \frac{R_{\Sigma}}{n}} K. 
\end{equation}
Such exponential rates have already been obtained in \cite{LuZhou2016,FeiChen2018} but under the corresponding stronger separation conditions~\eqref{LZ:condition} and~\eqref{cond:FC2018}. Nevertheless, the numerical constants in the exponential rate of  \cite{LuZhou2016} are tighter than ours.

In the high-dimensional setting $R_{\Sigma}>n$, the misclassification rate of relaxed $K$-means  decays at the slower exponential rate $\tfrac{\Delta^4 n}{\sigma^4 K R_{\Sigma}}< \frac{\Delta^2}{\sigma^2}$ when the distances between the means satisfy 
\begin{equation}\label{eq:condition_slower_rate}
 \sqrt{\frac{R_{\Sigma}}{n}} K \lesssim    \frac{\Delta^2}{\sigma^2}\lesssim  \frac{R_{\Sigma}}{n} K 
\end{equation}
Up to our knowledge, this moderate signal regime was not previously covered in the literature.
The discrepancy between the rates $\tfrac{\Delta^4 n}{\sigma^4 K R_{\Sigma}}$ and  $\frac{\Delta^2}{\sigma^2}$ may seem suboptimal. Yet, as we explain below, this discrepancy is inherent to the lack of knowledge of the location of the means of the clusters, and it seems unavoidable in our clustering setting.

Actually, let us consider the arguably simpler problem of Gaussian supervised classification with a two-class balanced partition, a common spherical covariance $\Sigma_{k}=\sigma I_p$ and opposite means $\mu_{-1}= -\mu_1$ uniformly distributed on the Euclidean sphere $\partial B(0,\Delta/2)$. More precisely, assume that we have $n$ labeled observations $(X_{a},Z_{a})\in \R^p\times\ac{-1,1}$, for $a=1,\ldots,n$ distributed as follows. The labels $Z_{1},\ldots,Z_{n}$ are i.i.d.\ with uniform distribution on $\ac{-1,1}$, a random vector $\mu\in\R^p$ is sampled  uniformly over the sphere $\partial B(0,\Delta/2)$ independently of  $Z_{1},\ldots,Z_{n}$, and, conditionally on $Z_{1},\ldots,Z_{n},\mu$, the $X_{a}$ are independent Gaussian random variables with mean $Z_{a}\mu$ and covariance $\sigma^2I_{p}$. Then, direct computations (see Appendix \ref{sec:bayes} for details) show that the classifier minimizing the probability of misclassification of a new observation is 
$$\widehat{h}(x)= \textrm{sign}\pa{\bigg\langle {1\over n}\sum_{a=1}^n Z_{a}X_{a},x\bigg\rangle}.$$
According to the invariance of the distribution by rotation, the probability of misclassification is given by
$$\P\cro{Z_{\rm new}\neq \widehat h(X_{\rm new})} =   \mathbb{P}\left[\Big\langle \frac{\mu}{\sigma}+\frac{\epsilon}{\sqrt{n}}, \frac{\mu}{\sigma}+ \epsilon'\Big\rangle <0  \right],$$
where $\mu=[\Delta/2,0,\ldots,0]\in \R^p$,  and $\epsilon,\epsilon'$ are two independent standard Gaussian vectors in $R^p$.
For $\tfrac{\Delta^2}{\sigma^2}\gtrsim (1\vee \frac{p}{n})$, this error is of exponential order $\Delta^2/\sigma^2$ whereas,  for $(1\vee \sqrt{\tfrac{p}{n}}) \lesssim \tfrac{\Delta^2}{\sigma^2}\lesssim (1\vee \tfrac{p}{n})$, the error is of exponential order $n\Delta^4/p\sigma^4$, see again Appendix \ref{sec:bayes} for details. Hence, even in this simpler toy model, the exponential rate with respect to $s^2={\Delta^2\over \sigma^2}\wedge {n\Delta^4\over 2 p \sigma^4}$ is intrinsic. To formalize this argument, one should prove rigorously a minimax lower bound for our clustering problem, but this is beyond the scope of this paper. 

When the common covariance $\Sigma$ is not spherical, we cannot argue anymore that the exponential rate $s^2$ is optimal. Actually, in the Gaussian supervised classification setting with known means ${\mu_{-1},\mu_{1}}$ and (common) covariance $\Sigma$, the probability of misclassification is known to decay exponentially fast with the Mahalanobis distance $d^2_{\Sigma}(\mu_{-1},\mu_{1})=(\mu_{1}-\mu_{-1})^T\Sigma^{-1}(\mu_{1}-\mu_{-1})$ rather than $\Delta^2/\sigma^2=\|\mu_{1}-\mu_{-1}\|^2/|\Sigma|_{op}$, see e.g. Section 9.5.1 in \cite{HDS}. 
Hence, we expect that the optimal rate of decay should involve $d^2_{\Sigma}$ instead of $\Delta^2/\sigma^2$, at least for Gaussian mixtures when the sample size is large. 
Since $K$-means criterion is tightly linked to isotropic (sub)Gaussian mixtures, it is unlikely that the relaxed $K$-means procedure enjoys an exponential decay with respect to $d_{\Sigma}^2$. In addition, we explain  in the discussion page
 \pageref{discussion:SNR} at the end of Section \ref{sec:outline1}, that the term  $m\Delta^4/\nu^4$ showing up in (\ref{SNR:subG}) is a variance term which  seems hard to avoid.
We emphasize yet that, despite these drawbacks, the exponential decay $e^{-c's^2}$ remains interesting in this setting, as it is "dimensionless", in the sense that it depends only on the effective rank $R_{\Sigma}$ of $\Sigma$ and not on the ambient dimension $p$.

  \medskip

\noindent {\bf Perfect recovery.}
Theorem \ref{thm1} ensures perfect recovery as soon as $s^2 \ag \log(n)\vee (n/m)$. This requirement exactly matches the requirement derived in Theorem 1 of \cite{MartinNIPS} and it is, to the best of our knowledge, the  sharpest known result for polynomial-time algorithms in high-dimension ($p\geq n$).  In the large sample size setting where $n\ag p^3K^2 \log(pK)$, \cite{VEMPALA2004} ensures  perfect recovery  for an iterative spectral clustering under the condition $\Delta^2 \ag \sigma^2\big(\sqrt{K\log(n)}+\log(n)\big)$ which is weaker than our $s^2 \ag \log(n)\vee (n/m)$  for  $K\gg \log(n)$.
Again, this lack of optimality is likely to be an artifact of the proof which is only valid when $s^2\ag K$. Actually, when $K\gg \log(n)$, the condition $s^2\ag K$ enforces  $e^{-c's^2} \ll 1/n$ .

\subsection{Unequal trace case}\label{sec:hetero}

As long as Condition \eqref{cond:Gamma} in Theorem \ref{thm1} is satisfied, the misclassification rates $err(\widehat G,G)$ is less than  $e^{-c's^2}$ as for the equal trace case. Let us first discuss some regimes and choice of $\widehat{\Gamma}$ under which \eqref{cond:Gamma} is valid, and then discuss the upper-bound.

\subsubsection{Choice of $\widehat \Gamma$}
First, observe that, when the covariance matrices $\Sigma_{k}$ have unequal traces, uncorrected convex $K$-means may not satisfy \eqref{cond:Gamma} if $\|\Gamma\|_V\geq  m\Delta^2/64$. Such a behavior is not an artifact of our proof techniques but is intrinsic for $K$-means as argued in Proposition 3 of \cite{MartinNIPS}.

When the covariance  matrices $\Sigma_{k}$ have unequal traces, we suggest to use in \eqref{pecok} the estimator $\widehat \Gamma$ introduced in \cite{pecok, MartinNIPS} and defined as follows.
For any $a,b\in \ac{1,\ldots,n}$, let 
$$V(a,b) := \max_{c,d \in \ac{1,\ldots,n}\setminus\{a,b\}}\Big|\langle X_{a}-X_{b}, \frac{X_{c}-X_{d}}{\|X_{c}-X_{d}\|_2}\rangle\Big|,$$ 
denote a measure of dissimilarity between $a$ and $b$. Then, let 
$\widehat{b}_1 = \argmin_{b \in \ac{1,\ldots,n}\setminus\{a\}}V(a,b)$ and $\widehat{b}_2:= \argmin_{b\in \ac{1,\ldots,n}\setminus\{a,\widehat{b}_1\}}V(a,b)$ denote the two indices most similar to $a$ with respect to this dissimilarity. Our estimator $\widehat{\Gamma}$ is defined by 
\begin{align}\label{estimator:Gamma}
\widehat{\Gamma}:= {\rm Diag}\pa{\langle X_{a}-X_{\widehat{b}_{1}}, X_{a}-X_{\widehat{b}_{2}}\rangle_{a=1,\ldots,n}}.
\end{align}

When $m>2$, denoting by $\gamma^2=L^2\max_{k}\tr(\Sigma_{k})$ the maximum scaled trace of the covariances  (with $L$ defined in Assumption A1),
Proposition 4 in \cite{MartinNIPS} ensures, that with probability higher than $1-c/n^2$
	\begin{align}
	|\widehat{\Gamma}-\Gamma|_\infty \al \left( \sigma^2 {\log n} + \sigma\gamma\sqrt{\log n} \right)\ . 
	\end{align}
Hence, Condition (\ref{cond:Gamma}) holds with probability larger than $1-c/n^2$ as soon as the condition
\begin{equation}
\label{eq:cond_gamma_simp1}
\Delta^2 \ag {\sigma^2\log(n)+\gamma\sigma \sqrt{\log(n)}\over m} 
\end{equation}
is met.
This condition can be compared to the condition $s^2\ag n/m$ arising in Theorem \ref{thm1}  by reformulating  it as
$${\Delta^2\over \sigma^2}\wedge {m\Delta^4\over \gamma^2\sigma^2} \ag {\log(n)\over m}.$$
In particular, we observe that this requirement is weaker than the condition $s^2\ag n/m$ when 
\begin{equation}\label{eq:cond_gamma_simp}
|\Sigma_{k}|_{op} \tr(\Sigma_{k}) \al \frac{n}{\log(n)} |\Sigma_{k}|_{F}^2\ ,\text{ for all }k=1,\ldots,K\ .  
\end{equation}
 This last condition is mild. For instance, it is met when
the ratio between the singular values $\sigma_{\ell}/\sigma_{1}$ of each $\Sigma_k$ decays faster than $1/\ell^{1+\epsilon}$ with $\epsilon>0$. It also holds when the condition number of each $\Sigma_{k}$ is upper bounded by $n/\log(n)$. 

\subsubsection{Discussion of the exponential rate}
Contrary to the equal trace case, even for spherical covariances $\Sigma_{k}=\sigma^2_{k}I_{p}$, we cannot argue anymore that our exponential rate $s^2$ is  optimal. Actually, in this case the signal-to-noise ratio $s^2$ should involve
$$\bar s^2=\min_{j\neq k} {\|\mu_{j}-\mu_{k}\|^2\over |\Sigma_{k}|_{op}\vee |\Sigma_{j}|_{op}}\quad\textrm{rather than}\quad {\min_{j\neq k} \|\mu_{j}-\mu_{k}\|^2\over \max_{k}|\Sigma_{k}|_{op}},$$
at least in the Gaussian case with large sample sizes. Some results on perfect recovery with respect to this  signal-to-noise ratio $\bar s^2$ have been derived in \cite{VEMPALA2004,AchlioptasMcSherry2005} in the asymptotic setting where $n\to \infty$. We do not know whether exponential decay with respect to $\bar s^2$ can be proved for the corrected relaxed $K$-means.

We emphasize that previously mentioned works~\cite{VEMPALA2004,AchlioptasMcSherry2005,mixon2016} are restricted to perfect recovery in the large sample size setting where (at least) $n\gg p$, while Theorem~\ref{thm1} provides a partial recovery bound that holds under  conditions \eqref{eq:cond_gamma_simp1} or \eqref{eq:cond_gamma_simp} which are non-asymptotic and can hold in any dimension $p$.
 As for~\cite{FeiChen2018}, they prove exponential decay with respect to $\tilde s^2={\min_{j\neq k} \|\mu_{j}-\mu_{k}\|^2/ \max_{k}|\Sigma_{k}|_{op}}$ in a higher signal regime and only for clusters of equal sizes.

\section{Clustering in Stochastic Block Models (SBM)}\label{sec:SBM}
\subsection{Stochastic Block Model (SBM)}\label{sec:SBM_definition}
We observe an undirected graph with $n$ nodes labeled by $a=1,\ldots,n$. We assume that the edges are independent and that  there exists an unknown partition $(G_{1},\ldots,G_{K})$ of the nodes and a symmetric matrix  $P\in [0,1]^{K\times K}$ such the probability to have an edge between $a\in G_{k}$ and $b\in G_{j}$ is equal to $P_{jk}$ when $a \neq b$ and 0 when $a=b$. In other words, the adjacency matrix $\X\in \R^{n\times n}$ is symmetric, with zero diagonal  and independent lower-diagonal entries fulfilling $\E[\X_{ab}]=P_{jk}$ for any $a\in G_{k}$ and $b\in G_{j}$, with $a\neq b$. 

Let us denote by $m_{k}$ the cardinality of $G_{k}$. 
There is a strong interest in the analysis of  recovery properties of SDP like
\begin{equation}\label{SDP:assortative}
\max_{B\in \mathcal{C}'}\langle \X,B\rangle,\quad \textrm{with}\ \mathcal{C}'=\ac{B: \textrm{PSD},\ B_{ab}\geq 0,\ B_{aa}=1,\ |B|_{1}=\sum_{k}m_{k}^2},
\end{equation}
or variants  of it. Such SDP are derived as a convex relaxation of the MLE optimization in the case where within group probability of connection $P_{kk}$ are all equal to $p$ and between group probability of connection $P_{jk}$ are all equal to $q$ with $q<p$, see e.g. \cite{chen2014statistical, 2015arXiv150705605P}. In the assortative setting, where the within group probabilities of connection are larger than the between group probabilities of connection, such SDP enjoy some very nice properties at all regimes, see  \cite{chen2014statistical,guedon2014community,2015arXiv150705605P,DAM15,2016arXiv160206410H,FeiChen2017}. 

Yet, we observe that SDP like (\ref{SDP:assortative}) seek for partitions maximizing within group connectivity. So they are  tied to the assortative setting and we cannot expect some good performances far away from this setting.

Instead of these SDP, we propose to solve  (\ref{relaxed-$K$-means}) or a variation of it, which is a relaxed-version of $K$-means applied to the adjacency matrix. The heuristic is as follows: two nodes $a$ and $b$ belonging to the same group share the same connectivity profile, that is, up to diagonal terms, the expectation of the columns $X_a$ and $X_b$ are equal. Therefore, it is tempting to recover the groups of the SBM by using distance clustering on the columns of the adjacency matrix $\X$.
In particular, (\ref{relaxed-$K$-means}) seeks for groups of nodes sharing similar connectivity profiles, instead of groups with maximal within group connectivity. The main difference between (\ref{relaxed-$K$-means}) and the SDP like (\ref{SDP:assortative}), is that the maximization is applied to $\X\X^T=\X^2$ instead of $\X$. So, compared to SDP like (\ref{SDP:assortative}), the SDP (\ref{relaxed-$K$-means}) seeks for a partition where the groups have a high density of common neighbors  rather than a high density of connections. Hence, the SDP (\ref{relaxed-$K$-means}) is not tied to the assortative case and can handle arbitrary matrices $P$.
We also point out, that, contrary to previous SDP investigated in the literature for SBM clustering, the relaxed $K$-means (\ref{relaxed-$K$-means}) does not require the knowledge  of the size of the groups, nor the knowledge of some parameters of  $P$.   
 Yet, in order to handle appropriately the sparse setting where $|P|_{\infty}=o(\log(n)/n)$, we need to add a constraint on $|B|_{\infty}$ in the program (\ref{relaxed-$K$-means}) to prevent the solution of $K$-means to produce too unbalanced partitions. As pointed in the previous section,  the norm $|B|_{\infty}$ of a matrix $B\in \mathcal{P}$ corresponding to a true partition is the inverse of the size of the smallest cluster. Thus, adding  a constraint on $|B|_{\infty}$ will avoid the formation of too unbalanced partitions.

Hence, we propose to solve the following constrained version of the relaxed $K$-means (\ref{relaxed-$K$-means})
\begin{equation}\label{constrained-relaxed-$K$-means}
\widehat B \in \argmax_{B\in \mathcal{C_{\alpha}}} \langle \X\X^T, B\rangle.
\end{equation}
where for $K/n\leq \alpha\leq 1$ 
$$\C_{\alpha}=\ac{B\in \R^{n \times n}:\ \textrm{Positive Semi Definite,}\ \tr(B)=K,\ B1=1,\ 0\leq B\leq \alpha}.$$
As explained in Theorem \ref{thm2} below, the parameter $\alpha$ can be chosen equal to 1 when $|P|_{\infty}\geq \log(n)/n$, but its choice is more constrained when $|P|_{\infty}\leq \log(n)/n$. We explain below Theorem~\ref{thm2} how $\alpha$ can be chosen in a data-driven way.

As for the mixture of sub-Gaussian, the  signal strength driving the exponential decay will be related to the Euclidean distance between two rows of the expected adjacency matrix. A row of $\X$ can be written as $X_{a}=\mu_{k}-P_{kk}e_{a}+E_{a}$, with 
$[\mu_{k}]_{b}=P_{jk}$  and $E_{a}=X_{a}-\E[X_{a}]$ for $a \in G_k$ and  $b \in G_{j}$. The square Euclidean distance between $\mu_{k}$ and $\mu_{j}$ is 
$$\|\mu_{k}-\mu_{j}\|^2=\sum_{\ell} m_{\ell} (P_{k\ell}-P_{j\ell})^2 \geq m \|P_{k:}-P_{j:}\|^2.$$
Similarly to the sub-Gaussian setting, we define 
\begin{equation}\label{Delta:SBM}
\Delta^2:= \min_{j\neq k}\Delta_{jk}^2\ ;\quad \Delta_{jk}^2:={\sum_{\ell} m_{\ell} (P_{k\ell}-P_{j\ell})^2}\geq {m} \|P_{k:}-P_{j:}\|^2
\end{equation}
which represents the signal strength in our analysis. 

We point out that $\|P_{k:}-P_{j:}\|\geq \sqrt{2}\,\lambda_{\min}(P)$, with $\lambda_{\min}(P)$ the smallest eigenvalue of $P$, so that $\Delta \geq \sqrt{2}\, m \lambda_{\min}(P)$.

Since the variance of a Bernoulli variable with small probability of success is roughly equal to this probability, we control the variance of $E_{a}$ with the following assumption.
\smallskip

\noindent {\bf Assumption A1'} We have $|P|_{\infty}\leq L$.
\smallskip

Under this assumption, we will prove that the misclassification error of the relaxed $K$-means decreases exponentially fast with the signal-to-noise ratio
$s^2=\Delta^2/L$.

\subsection{Main result}

We consider a partition $\widehat{G}$ obtained from 7-approximate $K$-medoid solution of $\widehat{B}$ as obtained in \eqref{constrained-relaxed-$K$-means}.
Next theorem provides an upper bound on the misclassified nodes decreasing exponentially fast with $s^2$. 

\begin{thm}\label{thm2}
Assume that Assumption A1' holds and set $s^2=\Delta^2/L$, with $\Delta$ defined by (\ref{Delta:SBM}).
Then, there exist three positive constants $c,c',c''$, such that for any
 $1/m\leq L \leq 1/\log(n)$,
\begin{equation}\label{eq:max-size}
{1\over m} \leq \alpha \leq \alpha(L):= {K^3\over n} e^{4nL}
\end{equation}
and
\begin{equation}\label{condition:SBM}
s^2 \geq c''n/m,
\end{equation} 
 with probability at least $1-c/n^2$,  the proportion of misclassified nodes is upper bounded by $$err(\widehat G,G)\leq e^{-c's^2}.$$
\end{thm}

We observe that we always have 
$$ \alpha(L)={K^3\over n} e^{4nL} \geq {K^3\over m}\times {m\over n}e^{4n/m}\geq {10^4\over m},$$ 
so the condition (\ref{eq:max-size}) is non-void. We discuss this condition into more details in the next section.

In practice, we can set $\alpha$ to the value $\widehat{\alpha}= {K^3\over n} e^{2nd_{{\X}}} \wedge 1$, where 
 $d_{{\X}}$ denotes the density of the graph. Next corollary provides a partial recovery bound when plugin   this data driven choice for  $\alpha$  in \eqref{constrained-relaxed-$K$-means}.
\begin{cor}\label{cor:SBM_alpha}
Assume that Assumption A1' holds, that $1/m\leq L \leq 1/\log(n)$, that the density $d_{\X}$ fulfills $\mathbb{E}[d_\X]\geq n^{-1}$ and
\begin{equation}\label{eq:max-size-2}
{1\over m} \leq {K^3\over n} e^{n\mathbb{E}[d_{\X}]}.
\end{equation}
Then, there exist three positive constants $c,c',c''$, such that when $s^2 \geq c''n/m$, and when setting $\alpha=\widehat \alpha$ in \eqref{constrained-relaxed-$K$-means},  then $err(\widehat G,G)\leq e^{-c's^2}$  with probability at least $1-c/n^2$.
\end{cor}
Indeed, $n(n-1)d_{{\X}}/2$ is stochastically dominated by a Binomial distribution with parameter $(n(n-1)/2, |P|_{\infty})$ so that $d_{{\X}}\leq 2|P|_{\infty}$ with probability higher than $1-e^{-0.3n(n-1)|P|_{\infty}}$ and hence $\widehat{\alpha}\leq \alpha(L)\wedge 1$. Conversely, Bernstein inequality together with $\mathbb{E}[d_\X]\geq n^{-1}$  ensures that $d_{{\X}}\geq \E[d_{\X}]/2$ with probability larger than $1-c/n^2$. Hence (\ref{eq:max-size-2}) ensures that $m^{-1}\leq \widehat \alpha$ with probability larger than $1-c/n^2$.

\subsection{Discussion} 
\noindent {\bf An SDP not tied to the assortative case.}
Our results cover a wide range of settings going beyond the assortative case usually handled by SDP algorithms. This is due to the fact that usual SDP criteria considered SBM clustering are derived as convex relaxation of MLE in the assortative case; while (\ref{relaxed-$K$-means}) is derived as a convex relaxation of $K$-means.
\medskip

\noindent {\bf On the condition (\ref{eq:max-size}).} 
The constraint $|B|_{\infty}\leq \alpha(L)$ is needed in our proofs in order to avoid the concentration of $\widehat B$ on the high-degree nodes. We observe first that $\alpha(L)\geq 1$ when 
$$L\geq {\log(n/K^3)\over 4n},$$ in which case we can take $\alpha=1$, which amounts to remove the constraint.

For smaller value of $L$ the constraint $|B|_{\infty}\leq \alpha(L)$ becomes active, with $\alpha(L)$ decreasing when $L$ decreases down to the extreme value $L=1/m$. We emphasize again that the condition on $\alpha$ does not require the knowledge of the true size of the groups but only constrains the size of the smallest group. Besides, the condition \eqref{eq:max-size} can be met as long as 
$ m \geq \tfrac{n}{K^3}e^{-4nL}$. For $L$ scaling as $l_0/n$, the size of smallest cluster can still be as small as $\tfrac{n}{K^3}e^{-4l_0}$ allowing for unbalanced partitions.

\medskip

\noindent {\bf Assortative case}. To start with, let us make explicit the value of $s^2=\Delta^2/L$ in the assortative case, with within group probabilities of connection $p$, between group probabilities of connection $q$ (with $q<p$) and balanced group sizes ($m\approx n/K$). In this case, $s^2=2m(p-q)^2/p$, and we obtain the same rate of exponential decay as in \cite{Chin15,YunP14b,AbbeSandon2015a,Gao2017, FeiChen2017}, but without the tight constants of ~\cite{YunP14b,Gao2017} in the exponential rate. We emphasize yet that we handle unknown group sizes, and with only a mild constraint on the group the sizes.
Besides, Theorem~\ref{thm2} ensures perfect recovery for
$${(p-q)^2\over p} \ag {K(K\vee \log(n))\over n},$$
matching the best known results (up to constants) for polynomial-time algorithms~\cite{chen2014statistical}.

\medskip

\noindent {\bf Partial recovery for General Model}. 
To the best of our knowledge, outside the assortative case, the only other exponentially decaying misclassification error is stated in Theorem 4 in \cite{AbbeSandon2015a} for a quite different procedure. Their results do not cover the same regime as ours, since they focus on the sparse regime where $P=P_0/n$ with $P_0$ a fixed matrix and $n\to \infty$. 
 For simplicity, let us discuss again the case of balanced groups where $m\approx n/K$. With our notation, Theorem~4 in \cite{AbbeSandon2015a} shows (under some conditions) that, in the sparse regime, the misclassification error is upper-bounded by
$e^{-c\tilde s^2}$ where
$$\tilde s^2= {m \lambda_{\min}(P)^2\over K\lambda_{\max}(P)}.$$
Since $\lambda^2_{\min}(P)\leq {2} \Delta^2/m$ and $\lambda_{\max}(P)\geq |P|_{\infty}$, we observe that
$$\tilde s^2\leq {2\Delta^2\over K |P|_{\infty}}={2\over K} s^2,$$
so that their  exponential decay with respect to $\tilde s^2$ is slower than the exponential decay with respect to $s^2$. Yet, this discrepancy between the rates is partly due to the fact that the exponential decay $e^{-c\tilde s^2}$ in \cite{AbbeSandon2015a} is valid  in regimes where our Condition \eqref{condition:SBM} for Theorem \ref{thm2} is not met. In other words, when the signal is low, \cite{AbbeSandon2015a} achieves the exponential decay $e^{-c\tilde s^2}$ which is not covered by our theory, while for stronger signals (where $s^2\ag K$ holds), our exponential decay $e^{-c s^2}$ is faster at least by a $K$ factor. We again emphasize that Theorem~\ref{thm2} is also valid in denser regime than that of \cite{AbbeSandon2015a}.

\medskip 

\noindent {\bf Perfect recovery for General Model}. 
From Theorem~\ref{thm2}, we derive that relaxed $K$-means achieves exact recovery as long as 
\begin{equation}\label{eq:exact_GSBM}
 s^2 \gtrsim \frac{n}{m}\vee \log(n)\ . 
\end{equation}

Again, the only other results we are aware of in the general model is from~\cite{AbbeSandon2015a} where the authors consider the asymptotic regime $P=\overline{P} \log(n)/n$ with $\overline{P}$ (and therefore also $K$) fixed and $n\to \infty$. In this setting,  they proved  that perfect recovery is possible in the balanced case ($m=K/n$)  if and only if
\begin{equation}\label{eq:abbe_sandon}
\lim_{n\to \infty} {m\over \log (n)} \min_{j\neq k} D_{+}(P_{j:}|| P_{k:}) >1\ , 
\end{equation}
where, for two vectors $q$ and $p$, 
\begin{equation}\label{eq:CH}
D_{+}(q||p) =\max_{t\in [0,1]} \sum_{x} p_{x} \pa{1-t+t {q_{x}\over p_{x}} -\pa{q_{x}\over p_{x}}^t}.
\end{equation}
Since $\Delta^2$ is based on the Euclidean distances between the columns of $P$ instead of $D_{+}$, our results cannot guaranty perfect recovery up to the exact threshold of \eqref{eq:abbe_sandon}. Yet, in the case where $\min_{j,k}P_{jk}/\max_{jk}P_{jk}$ is bounded away from zero, we can compare $\Delta^2$ to $\min_{j\neq k} D_{+}(P_{j:}|| P_{k:})$.
Actually, according to Lemma \ref{lem:CH} in Appendix \ref{section:CH}, we have
$$D_{+}(q||p) \leq {1\over 4\rho} \sum_{x} {(p_{x}-q_{x})^2\over p_{x}},\quad \textrm{when}\ \ \min_{x} {q_{x}\over p_{x}} \geq \rho >0.$$
Hence, when $\min_{j,k}P_{jk}\geq \rho \max_{j,k}P_{jk}$, we have
$$m\cdot \min_{j\neq k} D_{+}(P_{j:}|| P_{k:}) \leq {m\over 4 \rho^2} \min_{j\neq k} {\|P_{j:}-P_{k:}\|^2\over |P|_{\infty}} = {s^2\over 4\rho^2},$$
so that  $s^2$ and $m\cdot \min_{j\neq k} D_{+}(P_{j:}|| P_{k:})$  differs from at most a factor $1/(4\rho^2)$. As a consequence, in the asymptotic setting of~\cite{AbbeSandon2015a}, the condition \eqref{eq:exact_GSBM} achieves up to constants (depending on $\rho$) the optimal threshold~\eqref{eq:CH}.

When the ratio $\max_{jk}P_{jk}/\min_{jk} P_{jk}$ is unbounded, our SNR $s^2$ involving $|P|_{\infty}$ is no longer optimal. For example, in the weak assortative setting with clusters of equal sizes,  \cite{amini} derives some performance bounds which are  better  than ours when $\max_{jk}P_{jk}/\min_{jk} P_{jk}$ is unbounded.

\medskip

Finally, we point out that we do not need to de-bias the relaxed $K$-means as in (\ref{pecok}) for the sGMM. This is due to the fact that the size of the bias is small compared to the size of the fluctuations in this setting.

\section{Outline of the proofs}\label{sec:outline}

We write $B^*\in \mathcal{C}$ for the matrix associated to the true partition $G$ of the data set. Following the definition in Section \ref{sec:def:$K$-means}, we have  $B^*_{ab}=0$ unless $a$ and $b$ belong to the same group, in which case $B^*_{ab}=1/|G_k|$. 

\subsection{Outline of the proof of Theorem \ref{thm1}}\label{sec:outline1}
In this section, we describe the main lines of the proof of Theorem \ref{thm1}.  We refer to Section~\ref{sec:proof:thm1} for all the details. The proof relies on three main  arguments detailed below:
\begin{enumerate}
\item First, similarly to \cite{2015arXiv151208425C, FeiChen2017}, the misclassification proportion of the final $K$-medoid clustering (\ref{eq:$K$-medoid}) can be directly controlled by the $\ell^1$-norm $|\widehat B-B^*|_{1}$, see Section~\ref{sec:$K$-medoid}.
\item Second, by comparing $\langle \X\X^T-\widehat \Gamma,B^*\rangle \leq \langle \X\X^T-\widehat \Gamma,\widehat B\rangle$, we can upper-bound $|B^*\widehat B-B^*|_{1}$, (which is closely related to $|\widehat B-B^*|_{1}$) by some "noise" terms.
\item Third, a careful analysis of the noise terms provides the claimed result. Following \cite{FeiChen2017}, we use the key inequality $\sum_{i=1}^na_{i}b_{i}\leq \sum_{i=1}^{|b|_{1}} a_{i}$ for any $a_{1}\geq a_{2}\geq \ldots\geq a_{n}$ and $b_{1},\ldots,b_{n}\in [0,1]$, combined with tight upper-bounds on sums of ordered statistics. This bound involving ordered statistics is tighter than the classical $\sum_{i=1}^na_{i}b_{i}\leq |a|_{\infty} |b|_{1}$ combined with upper-bound on $\ell^{\infty}$-norms. We underline, that this classical reasoning based on $\ell^{\infty}$-norms  is not tight enough in order to handle partial recovery results in this setting.
\end{enumerate}

For simplicity, we assume throughout  this section that all the inter-cluster distances $\Delta_{jk}=\|\mu_{k}-\mu_{j}\|$ are equal to $\Delta=\min_{j\neq k}\Delta_{jk}$, for $j\neq k$.
We refer to Section~\ref{sec:proof:thm1} for  the general case.  As claimed in the first point above, our main task is to prove the bound  $|\widehat B-B^*|_{1}\leq ne^{-c's^2}$,  with probability at least $1-c/n^2$.

Let $A\in\R^{n\times K}$ denote  the membership matrix $A_{ik}=\1_{i\in G_{k}}$, let $\mathbf{\mu}\in \R^{K\times p}$ be the matrix whose $k$-th rows is given by $\mu_{k}$ and let $E\in \R^{n\times p}$ be the matrix whose $a$-th rows is given by $E_{a}$. The observed matrix $\X$ can then be written as $\X=A\mathbf{\mu} + E$. For any matrix $B\in \C$, expanding the product $\X\X^T$, we can decompose the scalar product
\begin{align*}
\langle \X\X^T-\widehat \Gamma,B^*-B\rangle &= \langle A\mathbf{\mu}\mathbf{\mu}^TA^T  ,B^*-B\rangle + \langle EE^T-\Gamma,B^*-B\rangle\\
&+\langle \Gamma-\widehat \Gamma,B^*-B\rangle+\langle A\mathbf{\mu}E^T+E\mathbf{\mu}^TA^T,B^*-B\rangle.
\end{align*}
The first term can be interpreted as a signal term, which is minimized in $\mathcal{C}$ at $B=B^*$. The three remaining terms involve the noise.
\smallskip

\noindent {\bf Signal term.}
Some basic algebra (see Lemma \ref{lem:signal} in Section~\ref{sec:proof:thm1}) shows that this writes as
\begin{equation}\label{eq:lower:signal}
\langle A\mathbf{\mu}\mathbf{\mu}^TA^T  ,B^*-\widehat B\rangle = {1\over 4} \Delta^2 |B^*-B^*\widehat B|_{1}.
\end{equation}
Since $B^*\in \C$, by definition of $\widehat B$ we have $\langle \X\X^T-\widehat \Gamma,B^*-\widehat B\rangle \leq 0$ 
and hence
$$ {1\over 4} \Delta^2 \,|B^*-B^*\widehat B|_{1} \leq \textrm{noise terms}\ .$$
With a suitable control of the three noise terms we can therefore hope to get a control on $|B^*-B^*\widehat B|_{1}$. Relying on the following inequality
\begin{equation}\label{eq:normes1}
|B^*-\widehat B|_{1}\leq {2n\over m}|B^*-B^*\widehat B|_{1}
\end{equation}
proved in Lemma~\ref{misc}, Section~\ref{sec:proof:thm1}, this will allows us in turn to control $|B^*-\widehat B|_{1}$. Let us explain how we can control each of the three noise terms, and how they contribute to the final result.
\smallskip

\noindent{\bf Quadratic terms.}
The quadratic term $\langle EE^T-\Gamma,\widehat B-B^*\rangle$ is the most delicate one. We observe first, that for any matrix $M$, the product $B^*M$ averages the rows of $M$ within the groups
$$[B^*M]_{ab}={1\over |G_{k}|} \sum_{c\in G_{k}} M_{cb},\quad \textrm{for all } a\in G_{k}\ .$$
As a consequence, the variance of the entries of  $B^*(EE^T-\Gamma)$ is reduced by a factor at least $\sqrt{m}$ compared to $EE^T-\Gamma$. So decomposing 
$$EE^T-\Gamma=B^*(EE^T-\Gamma)+(EE^T-\Gamma)B^*-B^*(EE^T-\Gamma)B^*+(I-B^*)(EE^T-\Gamma)(I-B^*),$$
we observe that we need to control three terms similar to $\langle B^*(EE^T-\Gamma),\widehat B-B^*\rangle$, and a last term $\langle (I-B^*) (EE^T-\Gamma)(I-B^*),\widehat B-B^*\rangle$. Since $B^*$ is a projection matrix, this last term involves the projection of $EE^T-\Gamma$ onto the orthogonal of the range of $B^*$. The three first terms benefits from the averaging effect described above, but not the last one. Instead, we control the last term as in \cite{pecok,MartinNIPS}
\begin{align}
\langle (I-B^*) (EE^T-\Gamma)(I-B^*),\widehat B-B^*\rangle &\leq |EE^T-\Gamma|_{op} |(I-B^*)(\widehat B-B^*)(I-B^*)|_{*} \nonumber \\
&\leq {1\over 2m} |EE^T-\Gamma|_{op} |B^*-B^*\widehat B|_{1}, \label{eq:nuclear}
\end{align}
see Lemma~\ref{misc} for the last inequality. Since $|EE^T-\Gamma|_{op}\al \nu^2 \sqrt{n} + \sigma^2 n$ with probability at least $1-c/n^2$, we obtain that (\ref{eq:nuclear}) is upper-bounded by $\Delta^2 \,|B^*-B^*\widehat B|_{1}/16$ under the assumption $s^2\ag n/m$. So this term is smaller than half of the signal term, and we can remove it at the price of losing a factor 2 in the signal level. We emphasize that the condition $s^2\ag n/m$ is exactly tailored to get this control, and hence it is fully driven by the upper-bound  $|EE^T-\Gamma|_{op}\al \nu^2 \sqrt{n} + \sigma^2 n$.

Let us now turn to the three terms of the form $\langle B^*(EE^T-\Gamma),\widehat B-B^*\rangle=\langle B^*(EE^T-\Gamma),B^*\widehat B-B^*\rangle$. The simple	 inequality $\langle A ,B\rangle \leq |A|_{\infty} |B|_{1}$ as in \cite{pecok, MartinNIPS} leads to the control 
\begin{equation}\label{martin1}
\langle B^*(EE^T-\Gamma),\widehat B-B^*\rangle \leq {\nu^2 \sqrt{\log(n)}+\sigma^2 \log(n)\over \sqrt{m}} |B^*-B^*\widehat B|_{1}\,,
\end{equation}
with high probability.
This control is
good enough to prove perfect recovery at the right scale, but it is too crude in order to exhibit  partial recovery rates. Instead, we adapt the clever analysis of \cite{FeiChen2017}, which relies on the upper-bound, 
$$\langle A,B\rangle \leq \sum_{j=1}^{|B|_{1}} A_{(j)},\quad \textrm{for any $B$ with } 0\leq B_{ab} \leq 1,$$
where $A_{(1)}\geq A_{(2)}\geq \ldots$ are the entries of $A$ ranked in decreasing order and 
$\sum_{j=1}^ba_{j}=a_{1}+\ldots+a_{[b]}+(b-[b])a_{[b]+1}$, where $[b]$ is the integer part of $b$.
Lemma~\ref{lem:order1} based on Hanson-Wright inequality provides a control of the sum of the ordered statistics of 
$ B^*(EE^T-\Gamma)$ ensuring that with probability at least $1-c/n^2$, the following inequality holds
\begin{align}
\lefteqn{\langle B^*(EE^T-\Gamma),B-B^*\rangle}\nonumber\\   &\al \pa{{\nu^2\over \sqrt{m}} {\sqrt{\log\pa{nK^3\over |B^*-B^*B|_{1}}}}\vee  \sigma^2\log\pa{nK^3\over |B^*-B^*B|_{1}}} |B^*-B^*B|_{1},\label{partial1}
\end{align}
simultaneously for all $B\in\C$ fulfilling $|B^*-B^*B|_{1}\al m$. 
The main difference compared to (\ref{martin1}), is that the $\log(n)$ has been replaced by something of the form $\log(n/ |B^*-B^*B|_{1})$. As it will appear clearly in the last step, moving from $\log(n)$ to $\log(n/ |B^*-B^*B|_{1})$ is the key to obtain the control $ne^{-c's^2}$ on $|B^*-B^*\widehat B|_{1}$. We point out that (\ref{partial1}) is only guaranteed for $B\in\C$ fulfilling $|B^*-B^*B|_{1}\al m$. 
So we need to first get such a bound. 

Contrary to  \cite{FeiChen2017}, we do not  use Grothendieck inequality for a preliminary control, but instead we apply a first time our upper-bound in order to get the rough bound $|B^*-B^*B|_{1}\leq nK^3e^{-c \sqrt{n/m}}\al m$ when $s^2\ag m/n$ and then apply again our analysis by using $|B^*-B^*B|_{1}\al m$.
We refer to Section~\ref{sec:proof:thm1:mainbound} below Lemma \ref{lem:cross} (Page \pageref{page:mainbound}) for the details.  
\smallskip

\noindent{\bf Gamma term.} The term $\langle \Gamma-\widehat \Gamma,\widehat B-B^*\rangle$ can be directly controlled by inequality (B12) in \cite{MartinNIPS}, which is recalled in Lemma~\ref{lem:Gamma} for convenience
$$\langle \Gamma-\widehat\Gamma,\widehat B-B^*\rangle \leq {2\over m} |\Gamma-\widehat \Gamma|_{V} |B^*-B^*\widehat B|_{1}.$$
The condition (\ref{cond:Gamma}) exactly ensures that the right hand side is upper-bounded by $\Delta^2 \,|B^*-B^*\widehat B|_{1}/32$. So this term is smaller than a quarter of the signal term, and we can again remove it at the price of loosing another factor 2 in the signal level.
\smallskip

\noindent{\bf Cross-products term.} It remains to upper-bound the cross-products term
$$\langle A\mathbf{\mu}E^T+E\mathbf{\mu}^TA^T,\widehat B-B^*\rangle =\sum_{j\neq k}\sum_{a\in G_{k},\ b\in G_{j}} \langle E_{a}-E_{b},\mu_{k}-\mu_{j}\rangle B_{ab}.$$
 We recall that we describe here the case where  the $\Delta_{jk}$ are all equal to $\Delta$ for $j\neq k$. The general case is treated in Section~\ref{sec:proof:thm1}. We observe first that $ \langle E_{a},\mu_{k}-\mu_{j}\rangle$ is sub-Gaussian SubG$(\sigma^2\Delta^2)$ and
$$ \sum_{j\neq k}\sum_{a\in G_{k},\ b\in G_{j}}  |B_{ab}| ={1\over 2} |B^*-B^*\widehat  B|_{1},$$
 see Lemma~\ref{misc} for this last equality.
Hence, building again on the inequality $\langle A,B\rangle \leq \sum_{j=1}^{|B|_{1}} A_{(j)}$ and deviation bounds for sub-Gaussian random variables, we 
obtain that
\begin{equation}\label{eq:cross-terms}
\langle A\mathbf{\mu}E^T+E\mathbf{\mu}^TA^T,\widehat B-B^*\rangle \leq  |B^*-B^*\widehat B|_{1} \sqrt{\sigma^2 \Delta^2 \log\pa{nK^3\over |B^*-B^*\widehat B|_{1}}},
\end{equation}
with probability at least $1-c/n^2$.
\smallskip

\noindent{\bf Conclusion.} Take for granted that we have $|B^*-B^*\widehat B|_{1}\al m$. Then combining the above bounds leads to 
$$ \Delta^2 |B^*-B^*\widehat B|_{1} \al \pa{{\nu^2\over \sqrt{m}} {\sqrt{\log\pa{nK^3\over |B^*-B^*\widehat B|_{1}}}}\vee  \sigma^2\log\pa{nK^3\over |B^*-B^*\widehat B|_{1}}} |B^*-B^*\widehat B|_{1},$$
with probability at least $1-c''/n^2$. This bound can be rewritten as 
$$s^2 \al {\log\pa{nK^3\over |B^*-B^*\widehat B|_{1}}}$$
and hence $|B^*-B^*\widehat B|_{1}\leq nK^3 e^{-c' s^2}$. 

In light of the three last bounds, we can traceback the terms contributing to the two parts of $s^2=(\Delta^2/\sigma^2)\wedge (m \Delta^4/\nu^4)$. \label{discussion:SNR} The second term $m \Delta^4/\nu^4$ in $s^2$ comes from the first term in the upper-bound (\ref{partial1}). The ratio $\nu^2/\sqrt{m}$ in  (\ref{partial1}) is driven by the variance term in Hanson-Wright inequality, so there is very little room for possible improvement. The first term $\Delta^2/\sigma^2$ in $s^2$ mainly comes from (\ref{eq:cross-terms}), and secondary from the second term in (\ref{partial1}). As above, the term $\Delta \sigma$ arising in (\ref{eq:cross-terms}) is a variance term, so there is again very little room for possible improvement on this side.

Recalling the inequality (\ref{eq:normes1}), the bound $|B^*-B^*\widehat B|_{1}\leq nK^3 e^{-c' s^2}$ obtained above ensures that 
$|B^*-\widehat B|_{1}\leq 2nK^3(n/m) e^{-c' s^2}$. This last inequality does not seem to meet our expectations. Yet, when $s^2 \ag n/m$, it enforces the targeted inequality $|B^*-\widehat B|_{1}\leq ne^{-c' s^2}$ for some smaller constant $c'$.

\subsection{Outline of the proof of Theorem~\ref{thm2}}\label{outline2}
Let us define $D$ the diagonal matrix with $D_{aa}=P_{kk}$ for $a\in G_{k}$ and $\X'=\X+D$.
We observe that $|\X\X^T-\X'(\X')^T|_{\infty}\leq |D|_{\infty} (|\X|_{\infty}+|\X'|_{\infty})\leq 2L$,
so, when $\Delta^2 \ag Ln/m$ we have 
$$|\langle \X\X^T-\X'(\X')^T, \widehat B-B^*\rangle|\leq 2L|B^*-\widehat B|_{1}\leq 4L{n\over m} |B^*-B^*\widehat B|_{1} \leq 0.05 \Delta^2 |B^*-B^*\widehat B|_{1}.$$
So we can replace $\X$ by $\X'$ in our analysis, since this term is smaller than a fraction of the signal term.

Since $\X'_{a}= \mu_{k}+E_{a}$, with 
$[\mu_{k}]_{b}=P_{jk}$  and $E_{ab}=X_{ab}-\E[X_{ab}]$ for $b \in G_{j}$, the proof of Theorem~\ref{thm2} follows similar lines as that of Theorem~\ref{thm1}. The main differences lies in the symmetry of $\X$ and the different stochastic control of the Bernoulli variables. We briefly sketch the main lines below. Again, for simplicity, we assume in this section that $\Delta_{jk}=\Delta$ for all $j\neq k$. 

As in the previous section, we have
$$\langle \X'(\X')^T, B^*-\widehat B\rangle = {1\over 4} \Delta^2 |B^*-B^*\widehat B|_{1} + \langle EE^T,B^*-\widehat B\rangle + \langle A\mu E^T + E \mu^TA^T ,B^*-\widehat B\rangle.$$

\noindent 
{\bf Cross-product terms.}
The cross-product terms are handled similarly as before, the main difference is that we rely on Bernstein inequality instead of sub-Gaussian deviations, producing an additional term. Since var$(E_{ai})\leq L$ and $|\mu|_{\infty}\leq L$, we get the same bound as before, with $\sigma^2$ replaced by $L$ and the additional term equal to $L$ times the sum of $|B^*-B^*\widehat B|_{1}$ ordered exponential random variables
\begin{align}\label{eq:cross-terms:sbm}
\langle A\mathbf{\mu}E^T+E\mathbf{\mu}^TA^T,\widehat B-B^*\rangle &\leq  |B^*-B^*\widehat B|_{1} \sqrt{L \Delta^2 \log\pa{nK^3\over |B^*-B^*\widehat B|_{1}}}\nonumber\\
&+ L  |B^*-B^*\widehat B|_{1} \log\pa{nK^3\over |B^*-B^*\widehat B|_{1}}.
\end{align}

\noindent{\bf Quadratic terms.} We use the same decomposition as before, except that $\Gamma$ is absent here. Actually the operator norm of $\Gamma$ is upper bounded by $nL$ which is smaller than the size of the fluctuations of $EE^T$ around $\Gamma$.

First, we consider the expression involving $(I-B^*)EE^T(I-B^*)$. When $L\geq \log(n)/n$, we can directly use the bound~(\ref{eq:nuclear}) since, in that regime, $|EE^T|_{op}\leq cnL$ (see e.g.~\cite{LeiRinaldo}). 
Unfortunately,  for a smaller $L$, the operator norm $|EE^T|_{op}$ can only be upper-bounded by $nL+\log(n)$ (up to a multiplicative constant). Actually, even for $L$ of order $1/n$, the nodes with the highest degree enforces an operator norm of size at least $\log(n)/\log \log(n)$. Being compelled to follow an alternative approach for bounding this expression, we use, classically for sparse graphs, a trimming argument which amounts to remove high degree nodes to the adjacency matrix. 
For the trimmed adjacency matrix we can now apply the bound~(\ref{eq:nuclear}) since the operator norm of this trimmed matrix is at most $cnL$. Then, it remains to upper bound the residual term by a $\ell^1/\ell^{\infty}$ bound. Relying on the box-constraint $|B|_{\infty}\leq \alpha(L)$, we will then guaranty that this residual remains under control.

As for the term $\langle B^*EE^T,\widehat B-B^*\rangle$, we need to control some quadratic forms of centered Bernoulli variables. We get a control of the right order by splitting them into pieces and considering apart different sub-cases. 
The symmetry of $\X$ induces some interlaced dependencies that must be handled with care. It is the main hurdle of the proof.

At the end of the day, we obtain a bound of the form 
\[\langle EE^T,B^*-\widehat B\rangle\leq 0.05 \Delta^2 |B^*-B^*\widehat B|_{1} + L  |B^*-B^*\widehat B|_{1} \log\pa{nK^3\over |B^*-B^*\widehat B|_{1}}\ ,\]
for the quadratic terms. Then, the end of the proof of Theorem \ref{thm2} follows the same line as those of Theorem \ref{thm1}.

\section{Proof of Theorem~\ref{thm1}}\label{sec:proof:thm1}

We provide in this section the full proof of Theorem~\ref{thm1}. We recall that $m=\min_k |G_k|$ stands for the size of the smallest group.
 In the sequel, $B_{aG_k}$ stands for $a$-th row of $B$ restricted to the column in $G_k$ whereas $B_{G_kG_l}$ stands for the restriction of $B$ to rows in $G_k$ and columns in $G_l$. 
For two indices $j$ and $k$, $(-1)^{j=k}$ equals $-1$ when $j=k$ and $1$ when $j\neq k$.

\subsection{A few useful formulas}
We start by gathering some useful formulas,  in particular relating $|B^*-B|_{1}$ to $|B^*-B^*B|_{1}$.
\begin{lem}\label{misc}
We have for any $B\in \mathcal{C}$
\begin{align*}
(B^*-B^*B)_{ab}&= {1\over m_{k}}-{1\over m_{k}}\sum_{c\in G_{k}} B_{cb} \quad\textrm{if}\ a,b\in G_{k},\\
&= -{1\over m_{k}}\sum_{c\in G_{k}} B_{cb} \quad\textrm{if}\ a\in G_{k}, b\notin G_{k}.
\end{align*}
and
$$|B^*-B^*BB^*|_{1}=|B^*-B^*B|_{1}=2\sum_{j \neq k} |B_{G_{j}G_{k}}|_{1},$$
and
$$|(I-B^*)B(I-B^*)|_{*}\leq {|B^*-B^*B|_{1}\over 2m}\ ,$$
and
$$ |B^*-B|_{1} \leq {2n\over m} |B^*-B^*B|_{1}.$$
Besides, for any $n\times n$  matrix $B$, we have 
\[
  |B^*-B^*B|_{1}\vee |B^*-BB^*|_{1}\leq |B^*-B|_{1}.
\]
\end{lem}
\noindent{\bf Proof of Lemma \ref{misc}.} The two first displays follows from direct computations. 
The third display is given by (57) in \cite{pecok}.
For the next to  last display, we observe that
$$|(I-B^*)B(I-B^*)|_{1}\leq n |(I-B^*)B(I-B^*)|_{*} \leq {n\over 2m} |B^*-B^*B|_{1}.$$
The claim follows from $B=B^*B+(I-B^*)B(I-B^*)+BB^*-B^*+B^*-B^*BB^*$ and hence
$$|B^*-B|_{1}\leq |B^*-B^*B|_{1}+|(I-B^*)B(I-B^*)|_{1}+|BB^*-B^*|_{1}+|B^*-B^*BB^*|_{1}.$$
For the last display, we use that for $a\in G_k$, $(B^*M)_{ab}= |G_k|^{-1} \sum_{c\in G_k}M_{cb}$. By triangular inequality, this implies that 
\[
 |B^*M|_{1}= \sum_{k=1}^K\sum_{a\in G_k}\sum_{b=1}^n \Big||G_k|^{-1} \sum_{c\in G_k}M_{cb}\Big|\leq |M|_1 \ . 
\]
Similarly $|MB^*|_{1} \leq |M|_{1}$ and the last display follows by taking $M=B^*-B$ and using $(B^*)^2=B^*$.
\ \hfill$\square$

The next two lemmas recall two useful probabilistic bounds.

\begin{lem}\label{HW} {\bf (Hanson-Wright inequality)}
Let $\eps$ be the vector obtained by concatenation of $\Sigma_{k(1)}^{-1/2}E_{1},\ldots,\Sigma_{k(n)}^{-1/2}E_{n}$. 

Under Assumption A1, the random vector  $\eps$ is sub-Gaussian SubG$(L^2I_{np})$ and for all $t>0$
$$\P\cro{\eps^TA\eps-\E[\eps^TA\eps]\geq L^2(|A|_{F}\sqrt{t}+|A|_{op}t)}\leq  e^{-ct}.$$
\end{lem}

We refer to \cite{rudelson2013hanson} for a proof of this lemma. Next lemma rephrases Lemma A1 in \cite{MartinNIPS}.

\begin{lem}\label{spectral}
Let $E$ be the $n\times p$ matrix $E^T=[E_{1},\ldots,E_{n}]$
Under assumption A1, we have for all $t>0$
$$\P\cro{\big|EE^T-\E[EE^T]\big|_{op} \geq \nu^2 \sqrt{t}+\sigma^{2}t}\leq 2 \times 9^n e^{-ct}.$$
\end{lem}

\subsection{Bounding $|B^*-B^*\widehat B|_{1}$}\label{sec:proof:thm1:mainbound}
We recall that $s^2=\min(\Delta^2/\sigma^2, m\Delta^4/\nu^4)$.
As explained in Section~\ref{sec:outline1}, the main step in the proof of Theorem~\ref{thm1} is to  prove that for $s^2 \ag {n\over m}$, we have with probability at least $1-c/n^2$
\begin{equation}\label{bound:delta}
|B^*-B^*\widehat B|_{1} \leq nK^3 \exp(-c s^2)\ .
\end{equation}
For any $B\in\C$, we have the decomposition $$\langle \X\X^T-\widehat \Gamma,B^*-B\rangle = \langle S+\Gamma-\widehat \Gamma+W+W',B^*-B\rangle$$
where for $a\in G_{k}$, $b\in G_{j}$,
\begin{align*}
S_{ab}&= -0.5 \|\mu_{k}-\mu_{j}\|^2\\
W'_{ab}&= \langle E_{a}-E_{b},\mu_{k}-\mu_{j}\rangle\\
W_{ab}&= \langle E_{a},E_{b}\rangle -\E[\langle E_{a},E_{b}\rangle].
\end{align*}

Since $B^*\in \C$, we have $\langle \X\X^T-\widehat \Gamma,B^*-\widehat B\rangle \leq 0$ by definition of $\widehat B$, and hence
\begin{equation}\label{eq:optim}
\langle S,B^*-\widehat B\rangle \leq \langle\Gamma-\widehat \Gamma,\widehat B-B^*\rangle+ \langle W,\widehat B-B^*\rangle+ \langle W',\widehat B-B^*\rangle.
\end{equation}

For the term in the left-hand side of (\ref{eq:optim}), direct computations combined with Lemma~\ref{misc} give the following evaluation of the signal.

\begin{lem}\label{lem:signal}
We set $b_{jk}=|B_{G_{k}G_{j}}|_{1}$ and  $\Delta_{jk}=\|\mu_{k}-\mu_{j}\|$. Then, for any $B\in\C$, we have
$$ \langle S,B^*-B\rangle= 0.5 \sum_{j\neq k} \Delta_{jk}^2 b_{jk}\geq {\Delta^2\over 4}\,|B^*-B^*B|_{1}.$$
\end{lem}

For the first term in the right-hand side of (\ref{eq:optim}), we lift from \cite{MartinNIPS} the next lemma (see (B12) in \cite{MartinNIPS}). Recall that $|D|_{V}=\max_{a}D_{aa}-\min_{a}D_{aa}$.
\begin{lem}\label{lem:Gamma}
For any diagonal $\widehat \Gamma$, and any $B\in \C$, we have
$$\langle \Gamma-\widehat\Gamma,B-B^*\rangle \leq {2\over m} |\Gamma-\widehat \Gamma|_{V} |B^*-B^*B|_{1}.$$
\end{lem}

For the second term  in the right-hand side of (\ref{eq:optim}), as explained in the Section~\ref{sec:outline1}, we decompose $W$ into 
$W=(I-B^*)W(I-B^*)+B^*W+WB^*-B^*WB^*$.

In order to control the scalar product involving $(I-B^*)W(I-B^*)$,
we get by combining Lemma~\ref{spectral} and Lemma~\ref{misc} the following bound.
\begin{lem} Under Assumption A1, we have with probability at least $1-c/n^2$
$$\langle (I-B^*)W(I-B^*), B-B^*\rangle \leq {\nu^2\sqrt{n}+\sigma^2 n \over m} |B^*-B^*B|_{1},$$
simultaneously for all $B\in\C$.
\end{lem}
Hence, when $s^2 \ag {n\over m}$ and when (\ref{cond:Gamma}) holds, we have $$\langle (I-B^*)W(I-B^*), \widehat B-B^*\rangle +\langle \Gamma-\widehat\Gamma,\widehat B-B^*\rangle \leq 0.75 \langle S,B^*-\widehat B\rangle.$$

The remaining terms involving $W$ are controlled by the next lemma proved in Section~\ref{sec:quadra}.

\begin{lem}\label{lem:quadratic}
We set $\delta=|B^*(B^*-B)|_{1}$ and we assume that Assumption A1  holds.
Then, with probability at least $1-c/n^2$, we have for all $B\in\C$
$$\langle B^* W,B- B^*\rangle \al {\delta\over \sqrt{m}} \pa{\nu^2\sqrt{\log(nK^3/\delta)}+\sigma^2\sqrt{\delta\vee 1}\ \log(nK^3/\delta)}.$$
The same bound holds for $\langle B^* WB^*, B-B^*\rangle$.
\end{lem}

It remains to control the last term in the right-hand side of (\ref{eq:optim}) with the following lemma, proved in Section~\ref{sec:cross}.

\begin{lem}\label{lem:cross}
We set $\delta=|B^*(B^*-B)|_{1}$.
Under Assumption A1, with probability at least $1-c/n^2$, we have for all $B\in\C$ 
\begin{align*}
\langle W',B-B^*\rangle &\al \sigma \sum_{j\neq k} \Delta_{jk} b_{jk} \sqrt{\log(nK/b_{jk})}\\
& \al \sigma \sqrt{ \langle S,B^*-B\rangle}  \sqrt{\delta\log(nK^3/\delta)}.
\end{align*}
\end{lem}

\noindent{\bf Conclusion.}\label{page:mainbound} Focusing on $\widehat{B}$, we  set $\delta=|B^*(B^*-\widehat B)|_{1}$.
Combining (\ref{eq:optim}) with 
$$\langle W,\widehat B-B^*\rangle = \langle (I-B^*)W(I-B^*),\widehat B-B^*\rangle + 2 \langle B^*W,\widehat B-B^*\rangle-\langle B^*WB^*,\widehat B-B^*\rangle,$$
and the five previous lemmas, we obtain that when $s^2 \ag {n\over m}$ and when (\ref{cond:Gamma}) holds, with probability at least $1-c/n^2$
 $$\langle S,B^*-\widehat B\rangle \al \sigma \sqrt{ \langle S,B^*-\widehat B\rangle}  \sqrt{\delta\log(nK^3/\delta)} \vee {\delta\over \sqrt{m}} \pa{\nu^2\sqrt{\log(nK^3/\delta)}+\sigma^2\sqrt{\delta\vee 1}\ \log(nK^3/\delta)}.$$
According to Lemma~\ref{lem:signal}, we have $\langle S,B^*-\widehat B\rangle\geq \Delta^2 \delta / 4$, and hence the previous bound ensures that
 \begin{eqnarray}
 {\Delta^2} &\al &  \sigma^2\log(nK^3/\delta) \vee{1\over \sqrt{m}} \pa{\nu^2\sqrt{\log(nK^3/\delta)}\vee  \sigma^2\sqrt{\delta\vee 1}\ \log(nK^3/\delta)} \nonumber \\ 
 &\al & \sigma^2\log(nK^3/\delta) \vee{1\over \sqrt{m}} \pa{\nu^2\sqrt{\log(nK^3/\delta)}\vee  \sigma^2\sqrt{\delta}\ \log(nK^3/\delta)}
 \label{borne1}
 \ .
 \end{eqnarray}
Since $|B^*|_{1}=|\widehat B|_{1}=n$ and $\delta\leq |B^*-\widehat B|_{1}$ (see Lemma~\ref{misc}), we always have $\delta\leq 2n$ by  triangular inequality. Let us prove that we actually have $\delta \al m$. From (\ref{borne1}) and $\delta\leq 2n$, we obtain
$$ {\Delta^2} \al  \sigma^2\sqrt{n\over m}\log(nK^3/\delta) \vee{1\over \sqrt{m}} \pa{\nu^2\sqrt{\log(nK^3/\delta)}}$$
and hence
$$\delta \leq nK^3 \exp\pa{-c\pa{\sqrt{m\over n}{\Delta^2\over \sigma^2}} \wedge {m \Delta^4\over \nu^4} }\leq nK^3 \exp\pa{-c \sqrt{m\over n}\ s^2}\leq nK^3 \exp(-c' \sqrt{n/m}),$$
where the last bound comes from $s^2 \ag m/n$. 
Hence, 
$$\delta\al m \,(n/m)^4  \exp(-c' \sqrt{n/m})\al m.$$

We can now conclude. Since $\delta \al m$, the bound (\ref{borne1}) gives 
$$ {\Delta^2} \al  \sigma^2\log(nK^3/\delta) \vee{1\over \sqrt{m}} \pa{\nu^2\sqrt{\log(nK^3/\delta)}}$$
from which follows 
$$\delta \leq nK^3 \exp(- c s^2) .$$
The proof of (\ref{bound:delta}) is complete.

\subsubsection{Proof of lemma \ref{lem:quadratic}}\label{sec:quadra}

In the following we use the notation $m_{k}=|G_{k}|$, and hence $m=\min_{k}m_{k}$. 
Since $B^*$ is a projection matrix, we observe that 
$\langle B^*W,B-B^*\rangle = \langle B^*W,B^*B-B^*\rangle.$

We have
$$\langle B^* W,B^*B-B^*\rangle=\sum_{k,j=1}^K \sum_{b\in G_{j}}z_{b}^{(j,k)}\beta_{kb}$$
where $\beta_{kb}=m_{k}|(B^*-B^*B)_{ab}|$ with $a\in G_{k}$ and 
$$z^{(j,k)}_b= {(-1)^{j=k}  \over m_{k}} \sum_{a \in G_{k}} W_{ab},\quad \textrm{for}\ b\in G_{j}\ .$$
We observe that $0\leq \beta_{kb} \leq 1$ and $|\beta_{kG_{j}}|_{1}=|(B^*-B^*B)_{G_{k}G_{j}}|_{1}=:b_{jk}$. 
Hence, writing $z^{(j,k)}_{(1)}\geq z^{(j,k)}_{(2)}\geq \ldots$ for the sequence $\ac{z^{(j,k)}_b: b\in G_{j}}$ ranked in decreasing order, we have
$$\langle B^* W,B^*B-B^*\rangle\leq \sum_{j,k=1}^K\sum_{u=1}^{b_{jk}} z^{(j,k)}_{(u)}\,,$$
with the convention that for $b=r+f$ with $r$ integer and $0\leq f < 1$,
\begin{equation}\label{eq:convention}
\sum_{u=1}^b a_{u}= (a_{1}+\ldots+a_{r})+ fa_{r+1}\leq { \pa{\sum_{u=1}^{r} a_{u}}\vee \pa{ \sum_{u=1}^{r+1} a_{u}}}\ . 
\end{equation}
and for $0\leq b < 1$, $\sum_{u=1}^b a_{u}\leq b a_{(1)}$.

We control the sum of ordered statistics by the next lemma proved at the end of this section.
\begin{lem}\label{lem:order1}
For any integer $q$ in $[1,m_{j}]$ and $t\geq 0$, we have 
$$\P\cro{\sum_{u=1}^q z_{(u)}^{(j,k)} \ag  \sqrt{q\over m} (\nu^2\sqrt{t}+\sigma^2t)}\leq C_{m_{j}}^q e^{-ct}.$$
\end{lem}
Let us choose $t_{q}= c'' q \log(nK/q)$. Since $C_{n}^q\leq (en/q)^q$, 
$$\sum_{q=1}^{m_{j}} C_{m_{j}}^q e^{-ct_{q}}\leq \sum_{q=1}^n e^{- c' q\log(nK/q)}\al {1\over (nK)^2},$$
where the last bound can be obtained e.g. by 
$$ \sum_{q=1}^n e^{- c' q\log(nK/q)}\leq  \sum_{q=1}^{\sqrt{n}} e^{- 0.5c' q\log(nK)}+ \sum_{q=\sqrt{n}}^n e^{- c' q\log(K)}\al {1\over (nK)^2}\ .$$

 Hence, with probability at least $1-c/(nK)^2$, we have simultaneously for all integers $1\leq q\leq m_{j}$,
$$ \sum_{u=1}^{q} z^{(j,k)}_{(u)} \al q  \sqrt{1\over m} \pa{\nu^2\sqrt{\log(nK/q)}+\sigma^2\sqrt{q}\log(nK/q)}\ .$$
From \eqref{eq:convention}, we deduce that 
$$ \sum_{u=1}^{b_{jk}} z^{(j,k)}_{(u)} \al b_{jk}  \sqrt{1\over m} \pa{\nu^2\sqrt{\log(\tfrac{nK}{b_{jk}})}+\sigma^2\sqrt{b_{jk}\vee 1}\log(\tfrac{nK}{b_{jk} \vee 1})}\ .$$
With a union bound over $j,k=1,\ldots,K$ we obtain that the inequality above holds simultaneously for all $j,k$ with probability at least $1-c/n^2$.
As a consequence, with Jensen inequality and $\sqrt{b_{jk}}\leq \sqrt{\delta}$, we get
\begin{align*}
\langle B^* W,B^*B-B^*\rangle&\al  \sum_{j,k=1}^K b_{jk}  \sqrt{1\over m} \pa{\nu^2\sqrt{\log(nK/b_{jk})}+\sigma^2\sqrt{\delta \vee 1}\log(nK/b_{jk})}\\
&\al  {\delta\over \sqrt{m}} \pa{\nu^2\sqrt{\log(nK^3/\delta)}+\sigma^2\sqrt{\delta \vee 1}\log(nK^3/\delta)}\ .
\end{align*}

The term $\langle B^*WB^*,B^*-B\rangle =\langle B^* WB^*,B^*(B^*-B)B^*\rangle$ can be handled in the same way as $\langle B^*W,B-B^*\rangle$, by noticing that
$|B^*(B^*-B)B^*|_{1}=|B^*(B^*-B)|_{1}$ according to Lemma~\ref{misc}.
 \bigskip

 \noindent 
{\bf Proof of Lemma \ref{lem:order1}:}
Let $\eps$ be the vector obtained by concatenation of $\Sigma_{k(1)}^{-1/2}E_{1},\ldots,\Sigma_{k(n)}^{-1/2}E_{n}$ which is sub-Gaussian SubG$(L^2I_{np})$. {Let $\otimes$ refer to the Kronecker product.}
For a subset $Q \subset G_{j}$ with cardinality $q$, we have 
$$\sum_{b\in Q} z^{(j,k)}_b ={(-1)^{j=k}\over m_{k}} \sum_{a\in G_{k},b\in Q} \pa{E_{a}^TE_{b}-\E[E_{a}^TE_{b}]}.$$
 {For a subset $Q\subset \{1,\ldots, n\}$, define $1_{Q}\in \{0,1\}^p$ such that $(1_{Q})_i=1$ if $i\in Q$.} Define 
$$A= {(-1)^{j=k}\over m_{k}}1_{G_{k}}1_{Q}^T$$
fulfilling $|A|_{F}=|A|_{op}= \sqrt{q/m_{k}}$. Then,  with $\Sigma^{(j,k)}=\Sigma_{k}^{1/2}\Sigma_{j}^{1/2}$, we have
$$\sum_{b \in Q} z^{(j,k)}_b= \eps^T(\Sigma^{(j,k)}\otimes A)\eps-\E\cro{\eps^T(\Sigma^{(j,k)}\otimes A)\eps}\ . $$
Hence, since $|\Sigma^{(j,k)}|_{F}\leq \nu^2/L^2$ and $|\Sigma^{(j,k)}|_{op}\leq \sigma^2/L^2$,  we have
$|\Sigma^{(j,k)}\otimes A|_{F}\leq \nu^2\sqrt{q}/(L^2\sqrt{m})$ and $|\Sigma^{(j,k)}\otimes A|_{op}\leq \sigma^2\sqrt{q}/(L^2\sqrt{m})$.
The bound then follows from Lemma~\ref{HW}, and the $C_{m_{j}}^q$ possible choices of subset $Q$.

\subsubsection{Proof of Lemma \ref{lem:cross}}\label{sec:cross}
To start with, we observe that 
$$\langle W',B-B^*\rangle =\sum_{j\neq k}\sum_{a\in G_{k},\ b\in G_{j}} \langle E_{a}-E_{b},\mu_{k}-\mu_{j}\rangle B_{ab}.$$
By symmetry, it is enough to control $\sum_{j\neq k}\sum_{a\in G_{k},\ b\in G_{j}} \langle E_{a},\mu_{k}-\mu_{j}\rangle B_{ab}$. 
Define $b_{jk}=|B_{G_{k}G_{j}}|_{1}$. Since $|B_{aG_{j}}|_{1}\leq 1$, we have 
\begin{align*}
\sum_{j\neq k}\sum_{a\in G_{k},\ b\in G_{j}} \langle E_{a},\mu_{k}-\mu_{j}\rangle B_{ab}&= 
\sum_{j\neq k}\sum_{a\in G_{k}} \langle E_{a},\mu_{k}-\mu_{j}\rangle |B_{aG_{j}}|_{1}\\
&\leq \sum_{j\neq k}\sum_{a=1}^{b_{jk}}z_{(a)}^{(j,k)},
\end{align*}
where  $z_{(1)}^{(j,k)}\geq z_{(2)}^{(j,k)} \geq\ldots$ corresponds to the values $(z_{a}^{(j,k)}=\langle E_{a},\mu_{k}-\mu_{j}\rangle: a\in G_{k})$ ranked in decreasing order.
Next lemma, proved at the end of this section, provides a control on the sum of the ordered statistics.

\begin{lem}\label{lem:order2}
Under Assumption A1,
with probability at least $1-c/n^2$, simultaneously for all $j\neq k$ and all integers  $1\leq q\leq m_{k}$
$$\sum_{a=1}^{q}z_{(a)}^{(j,k)}\al \sigma \Delta_{jk}q \sqrt{\log(nK/q)}\ .$$
\end{lem}

As a consequence, Cauchy-Schwarz and Jensen inequalities, and \eqref{eq:convention} ensure that
\begin{align*}
\langle W',B-B^*\rangle &\al \sigma \sum_{j\neq k} \Delta_{jk} b_{jk} \sqrt{\log(nK/b_{jk})}\\
&\al \sigma \sqrt{\sum_{j\neq k} \Delta_{jk}^2 b_{jk}} \sqrt{\delta\sum_{j\neq k} \delta^{-1} b_{jk}\log(nK/b_{jk})}\\
& \al \sigma \sqrt{ \langle S,B^*-B\rangle}  \sqrt{\delta \log(\sum_{j\neq k}nK/\delta)}\leq\sigma \sqrt{ \langle S,B^*-B\rangle}  \sqrt{\delta \log(nK^3/\delta)},
\end{align*}
which is the inequality claimed in Lemma \ref{lem:cross}.

It remains to prove Lemma \ref{lem:order2}. {Let $q\leq m_k$ denote a positive integer.}
Since the $(z^{(k,j)}_{a}:a\in G_{k})$ are independent and sub-Gaussian SubG$(\sigma^2\Delta_{jk}^2)$, 
for any $(a_{1},\ldots,a_{q})$ we have $\sum_{i=1}^q z^{(j,k)}_{a_{i}}$ sub-Gaussian SubG$(q\sigma^2\Delta_{jk}^2)$ under Assumption A1. Hence
$$\P\cro{\sum_{a=1}^{q}z^{(j,k)}_{(a)}>t} \leq C_{m_{k}}^{q} e^{-t^2/(2\sigma^2\Delta_{jk}^2q_{jk})}.$$
For $t_{q_{jk}}=c' \sigma \Delta_{jk} {q_{jk}}\sqrt{\log(nK/q_{jk})}$, we then have 
\begin{align*}
\P\cro{\exists q_{jk}\ :\ \sum_{a=1}^{q_{jk}}z^{(j,k)}_{(a)}>c' \sigma\Delta_{jk}q_{jk}\sqrt{\log(nK/q_{jk})}} &\leq  \sum_{j\neq k} \sum_{q_{j,k}=1}^{m_k} C_{m_{k}}^{q_{jk}} e^{-cq_{jk}\log(nK/q_{jk})}\\
&\al {K^2\over (nK)^c} \al {1\over n^2}
\end{align*}
for $c\geq 2$.

\subsection{Final clustering bound}\label{sec:$K$-medoid}
We recall that, similarly to \cite{FeiChen2017}, the final clustering is provided by a $7$-approximate $K$-medoids (\ref{eq:$K$-medoid}) on the rows of the matrix $\widehat B$ output by  (\ref{relaxed-$K$-means}) or (\ref{pecok}). Next lemma connects the misclassification error $err(\widehat G,G)$ defined by (\ref{misclassify:error}) to the $\ell^1$-norm $|B^*\widehat B-B^*|_{1}$.

\begin{lem}\label{lem:classif}
The proportion $err(\widehat G,G)$ 
 of misclassified points is upper bounded by
$$err(\widehat G,G) \leq 60 \pa{n\over m}^2 {|B^*\widehat B-B^*|_{1}\over n}.$$ 
\end{lem}

Combining Lemma~\ref{lem:classif} and (\ref{bound:delta}), we get that, when $s^2 \ag {n\over m}$, with probability at least $1-c/n^2$,
$$err(\widehat G,G) \leq 60 \pa{n\over m}^5 e^{-c' s^2}\leq e^{-c'' s^2}.$$ 
Theorem~\ref{thm1} then follows.

\subsubsection{Proof of Lemma \ref{lem:classif}}

The proof of Lemma \ref{lem:classif} is close to the proof of Proposition~3 in \cite{FeiChen2017}. We sketch below the main lines of this proof, referring to \cite{FeiChen2017} when the arguments are the same.

We define $\tilde B=\widehat A\widehat M$, with $(\widehat A,\widehat M)$ obtained in (\ref{eq:$K$-medoid}) and we define $A^{*}\in \R^{n\times k}$ by $A^*_{ak}=\1_{a\in G_{k}}$.
 We also define the sets
$T_{k}=\ac{a\in G_{k}: |B^*_{a:}-(\tilde BB^*)_{a:}|_{1}<1}$, $S_{k}=G_{k}\setminus T_{k}$, $R_{1}=\ac{k:T_{k}=\emptyset}$, 
$R_{2}=\ac{k: T_{k}\neq \emptyset\ \text{ and }\ \widehat A_{a:}=\widehat A_{b:}\ \forall a,b \in T_{k}}$ and 
$R_{3}=\{1,\ldots,K\}\setminus(R_{1}\cup R_{2})$.

The proof is decomposed into 4 steps.
\medskip

\noindent{\bf Step 1:} 
For $a\in T_{k}$ and $b\in T_{j}$ with $k\neq j$, we have $|B^*_{a:}-B^*_{b:}|_{1}\geq 2$ since $a\in G_{k}$ and $b\in G_{j}$ and $|B^*_{a:}-(\tilde B B^*)_{a:}|_{1}+ |B^*_{b:}-(\tilde B B^*)_{b:}|_{1}<2$ by definition of $T_{k}$ and $T_{j}$. Hence, 
$$|(\tilde B B^*)_{a:}-(\tilde B B^*)_{b:})|_{1} \geq |B^*_{a:}-B^*_{b:}|_{1} - |B^*_{a:}-(\tilde B B^*)_{a:}|_{1} - |B^*_{b:}-(\tilde B B^*)_{b:}|_{1}>0,$$
and so $\tilde B_{a:}\neq \tilde B_{b:}$ from which follows that $\widehat A_{a:}\neq \widehat A_{b:}$.

Hence, for $j,k\in R_{2}$, $a\in T_{k}$ and $b\in T_{j}$, we have $$\widehat A_{a:}\neq \widehat A_{b:}\quad \textrm{if and only if}\quad  j\neq k.$$ So all points in $\cup_{k\in R_{2}} T_{k}$ 
are well classified. Hence, there exists a permutation $\pi$ such that 
$$\left|\ac{a: \widehat A_{\pi(a):}\neq A^*_{a:}}\right| \leq n-\sum_{k\in R_{2}} |T_{k}|=S+ \sum_{k\in R_{3}}|T_{k}|\leq S+|R_{3}|n,$$
where $S=\sum_{k=1}^K|S_{k}|$. 
\medskip

\noindent{\bf Step 2:} The same arguments as in Claim 2 in \cite{FeiChen2017} ensures that $|R_{3}|\leq |R_{1}|$.
\medskip

\noindent{\bf Step 3:} We prove now  the inequalities $m|R_{1}|\leq S \leq |B^*-\tilde B B^*|_{1}$.
Actually, we have
$$S\geq \sum_{k\in R_{1}} |S_{k}|=\sum_{k\in R_{1}} |G_{k}|\geq m |R_{1}|$$
by definition of $R_{1}$ and $m=\min_{k}|G_{k}|$.
In addition,  since $|B^*_{a:}-(\tilde BB^*)_{a:}|_{1}\geq 1$ for $a\in \cup_{k}S_{k}$,
$$S\leq \sum_{k}\sum_{a\in S_{k}} |B^*_{a:}-(\tilde BB^*)_{a:}|_{1} \leq |B^*-\tilde BB^*|_{1}.$$

\noindent{\bf Step 4:}  The same arguments as in Claim 4 in \cite{FeiChen2017} ensures that $|B^*-\tilde B|_{1} \leq 15 |B^*-\widehat B|_{1}$.

We can now conclude. Combining the 3 first steps, we obtain that 
$$\left|\ac{a: \widehat A_{\pi(a):}\neq A^*_{a:}}\right|\leq S(1+n/m)\leq (1+n/m)  |B^*-\tilde BB^*|_{1}.$$
In addition, Step 4 and Lemma~\ref{misc} ensure that
$$|B^*-\tilde B B^*|_{1}\leq |B^*-\tilde B|_{1} \leq 15 |B^*-\widehat B|_{1} \leq 30{n\over m} |B^*-B^*\widehat B|_{1}.$$
The claim of Lemma \ref{lem:classif} then holds by combining the last two displays.

\section{Proof of Theorem \ref{thm2}} \label{sec:proof:thm2}
We provide in this section a full proof of Theorem~\ref{thm2}. 
The lines are very closed to those of the proof of Theorem~\ref{thm1}. In particular, all we need is to prove in (\ref{bound:delta})
in our setting. 
As explained in Section~\ref{outline2}, we have 
$$\langle \X\X^T, B^*-\widehat B\rangle - \langle \X'(\X')^T, B^*-\widehat B\rangle \leq 0.05 \Delta^2 |B^*-B^*\widehat B|_{1},$$
so 
$$\langle S,B^*-\widehat{B}\rangle= 0.5 \sum_{j\neq k} \Delta_{jk}^2 b_{jk}\leq  \langle EE^T,\widehat B-B^*\rangle+ \langle W',\widehat B-B^*\rangle+0.05 \Delta^2 |B^*-B^*\widehat B|_{1},
$$
with $W'_{ab}= \langle E_{a}-E_{b},\mu_{k}-\mu_{j}\rangle$ for $a\in G_{k}$ and $b\in G_{j}$.

The cross-product $\langle W',\widehat B-B^*\rangle$ can be easily bounded with a variation of Lemma~\ref{lem:cross}. 
\begin{lem}\label{lem:cross:sbm}
With probability at least $1-c/n^2$, we have for all $B\in\C$
\begin{align*}
\langle W',B-B^*\rangle 
& \al  \sqrt{  \langle S,B^*-B\rangle}  \sqrt{\delta L \log(nK^3/\delta)}+ L  \delta \log\pa{nK^3\over \delta}\ , 
\end{align*}
 {where $\delta=|B^*(B^*-B)|_{1}$.}
\end{lem}

The term $ \langle EE^T,\widehat B-B^*\rangle$ must be handled with more care.
We first focus on $ \langle (I-B^*)EE^T(I-B^*),\widehat B-B^*\rangle$.
\begin{lem}\label{lem:trimmed}
When $\Delta^2\ag Ln/m$, with probability at least $1-c/n^2$, we have 
$$\langle (I-B^*)EE^T(I-B^*),\widehat B-B^*\rangle \leq  0.05 \Delta^2 |B^*-B^*\widehat B|_{1},$$
when $L\geq \log(n)/n$ and 
$$\langle (I-B^*)EE^T(I-B^*),\widehat B-B^*\rangle \leq  0.05 \Delta^2 |B^*-B^*\widehat B|_{1}+0.04 n (nL)^2 e^{-cnL},$$
when $L\leq \log(n)/n$ and $\alpha(L)\leq K^3n^{-1} e^{cnL}$.
In particular, either $|B^*-B^*\widehat B|_{1}\leq nK^3 e^{-c \Delta^2/4L}$ or 
$$\langle (I-B^*)EE^T(I-B^*),\widehat B-B^*\rangle \leq  0.06 \Delta^2 |B^*-B^*\widehat B|_{1}.$$
\end{lem}
It remains to control the average terms of the form $ \langle B^*EE^T,\widehat B-B^*\rangle$.

\begin{lem}\label{lem:horrible}
When $\Delta^2\ag Ln/m$, with probability at least $1-c/n^2$, we have 
$$\langle B^*EE^T,\widehat B-B^*\rangle \leq 0.05 \Delta^2 \delta + c' L  \delta \log\pa{nK^3\over \delta}\ .$$
\end{lem}

Putting the three last lemmas together, we conclude that, either $\delta=|B^*-B^*\widehat B|_{1}$ fulfills $\delta \leq nK^3 e^{-c\Delta^2/4L}$, or with probability larger than $1-c''/n^2$, we have
$$ \delta \Delta^2 \al \delta L \log (nK^3/\delta).$$
In any case, it follows that $\delta \leq nK^3 e^{-c' \Delta^2/L}$ with probability at least $1-c/n^2$, which gives (\ref{bound:delta}). We conclude the proof of Theorem \ref{thm2} by following the same lines as for Theorem \ref{thm1}.

\subsection{Proof of Lemma \ref{lem:cross:sbm}}
As in the proof of Lemma \ref{lem:cross}, we denote by $z_{(1)}^{(j,k)}\geq z_{(2)}^{(j,k)} \geq\ldots$  the values $(z_{a}^{(j,k)}=\langle E_{a},\mu_{k}-\mu_{j}\rangle: a\in G_{k})$ ranked in decreasing order and we have
$$\langle W',B-B^*\rangle \leq \sum_{j\neq k}\sum_{a=1}^{|B_{G_{j}G_{k}}|_{1}}z_{(a)}^{(j,k)}\ ,$$
{where we use  the same convention as in \eqref{eq:convention}  when  $|B_{G_{j}G_{k}}|_{1}$ is not an integer.}
For any $Q\subset G_{k}$ with cardinality $q$ we have
$$\sum_{a\in Q}z_{a}^{(j,k)}= \sum_{a \in Q} \sum_{i=1}^n \cro{(E_{ai}\1_{i<a}+E_{ia}\1_{i>a})(\mu_{k}-\mu_{j})_{i}}.$$
Since each variable $E_{ab}$ for $a<b$ appears at most twice, $\sum_{i}$var$(E_{ai}(\mu_{k}-\mu_{j})_{i})\leq L \Delta_{jk}^2$ and $|E_{ai}(\mu_{k}-\mu_{j})_{i}|\leq 2L$, Bernstein inequality ensures that 
$$\P\cro{\sum_{a\in Q}z_{a}^{(j,k)}\ag \sqrt{qL\Delta_{jk}^2t}+Lt} \leq e^{-t}.$$
The conclusion of Lemma \ref{lem:cross:sbm} is then derived by following the same lines as in the proof of Lemma \ref{lem:cross}.

\subsection{Proof of Lemma \ref{lem:trimmed}}

We first consider the case  $L\geq \log(n)/n$ and then turn to the sparse case $L\leq \log(n)/n$. When $L\geq \log(n)/n$, according e.g.  to
Theorem~5.2 in \cite{LeiRinaldo} we have $|EE^T|_{op}\al nL$, so (\ref{eq:nuclear}) ensures that when $\Delta^2\ag Ln/m$
$$ \langle (I-B^*)EE^T(I-B^*),\widehat B-B^*\rangle\al {nL\over m} |B^*-B^*\widehat B|_{1} \leq 0.05 \Delta^2 |B^*-B^*\widehat B|_{1}.$$
Let us now consider the case $L\leq \log(n)/n$.
As mentioned in Section~\ref{outline2}, we use a trimming argument. Let $\X^{tr}$ be the matrix $\X$ where we have removed the nodes with degrees larger than $\gamma nL$, with $\gamma=2^8+2$,  and set $E^{tr}=\X^{tr}-\E[\X]$. It is known that removing the high-degree nodes drastically reduces the operator norm of the adjacency matrix. For instance,
 Lemma~5 in \cite{FeiChen2017} ensures that $|E^{tr}(E^{tr})^T|_{op}\al nL$ so when $\Delta^2\ag Ln/m$, we have
\begin{equation}\label{eq:1_0}
 \langle (I-B^*)E^{tr}(E^{tr})^T(I-B^*),\widehat B-B^*\rangle \leq 0.05 \Delta^2 |B^*-B^*\widehat B|_{1}\ . 
\end{equation}

Similarly as in \cite{FeiChen2017}, we bound now the residual terms with a $\ell^1/\ell^{\infty}$ bound. Compared to  \cite{FeiChen2017}, the main additional difficulty comes from the  quadratic residuals, whereas the SDP in \cite{FeiChen2017} was only linear in $\X$. 
The first step is to bound the $\ell^1$ norm of the residual terms by the sum of the square degrees of the trimmed nodes. We start from 
\begin{eqnarray}
 \lefteqn{ \langle EE^T-E^{tr}(E^{tr})^T, (I-B^*)\widehat B(I-B^*)\rangle}&& \nonumber\\ &\leq & |(I-B^*)\widehat{B}(I-B^*)|_{\infty}|EE^T-E^{tr}(E^{tr})^T|_{1}  \nonumber \\ 
  &\leq & 2\alpha(L)|EE^T-E^{tr}(E^{tr})^T|_{1}   \nonumber \\
  &\leq & 2\alpha(L)\left(|(E-E^{tr})^2|_1 + 2 |(E-E^{tr})E^{tr}|_1 \right) \label{eq:upper_noise_sparse}
\end{eqnarray}
since $|\widehat{B}|_{\infty}\leq \alpha(L)$. 
\medskip 

\noindent {\bf Control of $|(E-E^{tr})^2|_1$.} The matrix $E-E^{tr}=\X-\X^{tr}$ is the adjacency matrix of the graph 
where we have only kept the edges involving at least one node with with degree larger than $\gamma nL$ and their neighbors. The $\ell^1$ norm $|(E-E^{tr})^2|_1$ then counts the number of paths of size $2$ in this graph. Write $\cT$ for the set of nodes with degree larger than $\gamma nL$. To evaluate the number of paths $(i_{1},i_{2},i_{3})$ of size $2$ in this graph, we consider apart the two cases $i_{2}\in \cT$ and $i_{2}\notin \cT$.

 Consider first the case  where the node $i_2$ belongs to $\cT$. Since both $i_1$ and $i_3$ are neighbors of $i_2$, we have at most $\sum_{i\in \cT} d^2(i)$ such paths, where $d(i)$ is the degree of $i$. Consider now the case where $i_2\notin \cT$. In this case, $i_1$ and $i_{3}$ belong to $\cT$. Since the degree of $i_2$ is less than $\gamma nL\leq d(i_{1})$, the number of such paths is  again smaller than $\sum_{i\in \cT} d^2(i)$. So, we have the bound
\[
 |(E-E^{tr})^2|_1 \leq 2\sum_{i\in \cT} d^2(i).
\]

\noindent {\bf Control of $|(E-E^{tr})E^{tr}|_1$.}
We have 
$$|(E-E^{tr})E^{tr}|_1\leq |(E-E^{tr})\X^{tr}|_1+|(E-E^{tr})\E[\X]|_1,$$
and we bound separately the two terms in the right hand side of the above inequality. 
First, we notice that $|(E-E^{tr})\X^{tr}|_1$ corresponds to the number of size $2$ paths $(i_1,i_2,i_3)$ such that $i_1$ belongs to $\cT$ and both $i_2$ and $i_3$ do not belong to $\cT$. Since $d(i_2)< d(i_1)$, we have again  $|(E-E^{tr})\X^{tr}|_1\leq \sum_{i\in \cT} d^2(i)$. 
As for $|(E-E^{tr})\E[\X]|_1$, this corresponds to the sum of the weights $(E-E^{tr})_{i_{1}i_{2}}\E[\X_{i_{2}i_{3}}]$ associated to paths $(i_1,i_2,i_3)$. Since the weight $\E[\X_{i_{2}i_{3}}]$ of the edge $(i_2,i_3)$ is less than $L$, the total weight of paths starting from $i_1\in \cT$ is upper-bounded by $nL\sum_{i\in \cT}d(i)$ and the total weight of paths starting from $i_1\notin \cT$ is also upper-bounded by $nL\sum_{i\in \cT}d(i)$, so $|(E-E^{tr})\E[\X]|_1\leq 2nL\sum_{i\in \cT}d(i)\leq  \sum_{i\in \cT}d^2(i)$ since $\cT$ is made of high-degree nodes. 

Coming back to \eqref{eq:upper_noise_sparse}, we conclude that 
\begin{equation}\label{eq:control_noise_B*}
  \langle EE^T-E^{tr}(E^{tr})^T, (I-B^*)\widehat B(I-B^*)\rangle \leq  8\alpha(L) \sum_{i\in \cT}d^2(i).
\end{equation}
It remains to bound the sum of the squared highest degrees. 
\medskip

\noindent {\bf Control of  $ \sum_{i\in \cT}d^2(i)$. }
We control the sum with a stratification argument. First, we shall get rid of the dependencies in $\X$ that are due to the symmetry. 
$$d^2(i)=\pa{\sum_{j: j>i}\X_{ij}+\sum_{j:j<i}\X_{ij}}^2\leq 2\pa{\sum_{j: j>i}\X_{ij}}^2+2\pa{\sum_{j:j<i}\X_{ij}}^2.$$
For a node $i$, we write $d_1(i)= \sum_{j: j>i}\X_{ij}$ and $d_2(i)= \sum_{j: j<i}\X_{ij}$. As a consequence, 
\[
 \sum_{i\in \cT}d^2(i)\leq 4\sum_{i=1}^n d^2_1(i) \1_{d_1(i)\geq \gamma Ln/2}+  4\sum_{i=1}^n d^2_2(i) \1_{d_2(i)\geq \gamma Ln/2}
\]
We focus on the first term, the second term can be bounded in the same way by symmetry. The following technical Lemma is stated in general form as it will be applied several times in the manuscript. Henceforth, $\log_2$ refers to the binary logarithm. 

\begin{lem}\label{lem:peel}
Consider any $\ell>0$ such that $\ell L\geq 1$.
Let $I\subset \ac{1,\ldots,n}$,  $J_{i}\subset \ac{i+1,\ldots,n}$, with $|J_{i}|\leq \ell$ and $S_{i}=\sum_{j\in J_{i}} E_{ij}$ for $i=1,\ldots,n$.
For $r_{0}\geq 2$ and any integer $r\geq 1$ we set $y_{r}=2^{r_{0}+r} \ell L$ and
$I_{r}=\ac{i\in I: y_{r} < S_{i} \leq y_{r+1}}$.
Then, for $1/14\leq \tau\leq 1/2$, we have
$$\P\cro{\bigcap_{r\geq 1}\ac{ |I_{r}|\leq 2n 2^{-\tau (r+r_{0}-2)y_{r}}}}\geq  1- {1+ \log_{2}(\tau^{-1}\log(2n))\over (2n)^{(1-\tau)/4\tau}}\ ,$$
for $n\geq 2$.
In addition, we always have $S_{i}\geq -\ell L$ and for $I_{0}=\ac{i\in I: S_{i}\leq y_{1}}$
$$\P\cro{\sum_{i\in I_{0}} S_{i}^2 \al n\ell L} \geq 1-1/n^3.$$
 \end{lem}

Since $\X_{ij}\in [0,1]$,  $d_1(i)\leq nL+ \sum_{j\geq i}E_{ij}$. Take $r_0=6$, $I=\{1,\ldots, n-1\}$ and  $J_i=\{i+1,\ldots, n\}$ for $i\in I$. Since we restrict ourselves to indices $i$ such that $d_1(i)\geq \gamma nL/2$, our choice of $\gamma$ implies that $S_i> y_1$. Taking $\ell=n$ and $\tau= 1/10$ in Lemma \ref{lem:peel},
 we obtain with probability at least $1-c/n^2$ that 
\begin{align*}
\sum_{i=1}^n d^2_1(i)\1_{d_1(i)\geq \gamma nL/2}& \leq \sum_{r=1}^{\infty}|I_r| (nL+2^{r_{0}+r+1} nL)^2\\ 
&\leq 2^8 n(nL)^2 \sum_{r=1}^{\infty} 2^{-(r_{0}+r-2) (nL 2^{r_{0}+r}/10-2)}\\
& \leq 2^8 n(nL)^2 \sum_{r=1}^{\infty} 2^{-11.8(r_{0}+r-2) nL }\\
&\leq 2^9 n(nL)^2 e^{-40nL}. 
\end{align*}

So putting pieces together with \eqref{eq:upper_noise_sparse}, we obtain 
$$ \langle EE^T-E^{tr}(E^{tr})^T, (I-B^*)\widehat B(I-B^*)\rangle \leq 2^{12} \alpha(L) n(Ln)^2 e^{-40nL} \leq 0.04(Ln)^2K^3 e^{-4nL },$$
where the last inequality holds when $L\geq 1/m\geq 2/n$ and
$ \alpha (L) \leq \tfrac{K^3}{n}e^{4nL}$. Then, with \eqref{eq:1_0}, we conclude that 
\begin{equation}\label{eq:00}
 \langle (I-B^*)EE^T(I-B^*), \widehat B-B^*\rangle \leq 0.05 \Delta^2 |B^*-B^*\widehat B|_{1}+  0.04(Ln)^2K^3 e^{-4nL }
\end{equation}

Let us prove the last statement of the lemma. Assume that $\delta=|B^*-B^*\widehat B|_{1}\geq nK^3 e^{-\Delta^2/L}$ and hence
$$\delta \Delta^2 \geq nLK^3 (\Delta^2/L)e^{- \Delta^2/L}.$$
Since $n/m\al \Delta^2/L \leq 4nL$ and $xe^{-x}$ is decreasing for $x>1$, then 
$(\Delta^2/L)e^{- \Delta^2/L}\geq 4nL e^{-4nL}$
and  $\delta \Delta^2 \geq 4 K^3 (nL)^2e^{-4nL}$. Coming back to \eqref{eq:00} concludes the proof.

\subsection{Proof of Lemma~\ref{lem:horrible}}
In order to properly handle the dependences between the symmetric entries of $E$ we split $E$ into two parts $E=U+U^T$ where the upper triangular matrix $U$ is such that $U_{ab}=E_{ab}$ for $a<b$ and $U_{ab}=0$ else. We have $E^2=U^2+(U^T)^2+UU^T+U^TU$ and by symmetry we only need to control 
$\langle B^*U^TU,\widehat B-B^*\rangle$ and $\langle B^*(U^T)^2,\widehat B-B^*\rangle$.
\medskip

\subsubsection{Case $U^TU$}
As in the proof of Lemma~\ref{lem:quadratic}, all we need is to prove the following bound. For $1\leq j,k\leq K$ and $b\in G_j$, define 
\[
z^{(j,k)}_b=  \ac{{\pa{-1}^{j=k}\over m_{k}} \sum_{i=1}^n\sum_{a\in G_{k}} U_{ia}U_{ib}:b\in G_{j}}\ .
\]
Let $z^{(j,k)}_{(1)}\geq z^{(j,k)}_{(2)}\geq \ldots$ be the random variables $z^{(j,k)}_b$ ranked in  decreasing order.
\begin{lem}\label{lem:quadratic:sbm}
There exists an event $\Omega$ of probability at least $1-c/n^2$, such that for all $1\leq j,k\leq K$ and  for any integer $q \in [1,m_{j}]$ and $t\geq 0$, we have 
 \begin{equation}\label{sum:ordered:sbm}
 \P\cro{\Omega\cap\ac{\sum_{u=1}^q z_{(u)}^{(j,k)}\ag L\pa{{q n\over m_k}+t}}}\leq 3C_{m_{j}}^qe^{-c't}.
 \end{equation}
\end{lem}

To conclude from Lemma~\ref{lem:quadratic:sbm}, we simply apply as in Lemma~\ref{lem:quadratic} a union bound
\begin{align*}
\P\cro{\exists q:\sum_{u=1}^qz^{(j,k)}_{(u)}\ag {nLq\over m_{k}}+Lq\log(\frac{nK}{q})}
&\leq \P[\Omega^c]+\P\cro{\Omega \cap \ac{\exists q:\sum_{u=1}^qz^{(j,k)}_{(u)}\ag {nLq\over m_{k}}+Lq\log(\frac{nK}{q})}}\\
&\leq  {c\over n^2}+3\sum_{j,k} \sum_{q_{j,k}=1}^n C_{m_{k}}^{q_{jk}} e^{-cq_{jk}\log(\frac{nK}{q_{jk}})}\\
&\al {1\over n^2}+{K^2\over (nK)^c} \al {1\over n^2}\ ,
\end{align*}
for $c\geq 2$. 
We denote $b_{jk}=|B_{G_{k}G_{j}}|_{1}$ so that $\delta=\sum_{j,k} b_{jk}$.
Then, arguing as in Lemma~\ref{lem:quadratic}, we use order variables with the convention~\eqref{eq:convention} together with Jensen inequality to conclude
\begin{eqnarray*}
\langle B^*UU^T,\widehat B-B^*\rangle &\leq& \sum_{k,j=1}^K \sum_{u=1}^{b_{jk}} z_{(u)}^{(j,k)}\\
&\al &\sum_{k,j=1}^K b_{jk}L \left(\frac{n}{m_{k}}+\log(\frac{nK}{b_{jk}})\right) \\
&\al & \delta L \left(\frac{n}{m}+L\log(\frac{nK^3}{\delta})\right)
\end{eqnarray*}
Since we assume that $\Delta\gtrsim Ln/m$, we have proved the desired bound. 

\medskip

\noindent{\bf Proof of Lemma~\ref{lem:quadratic:sbm}.}
With the notation of Lemma~\ref{lem:peel}, for $1\leq k\leq K$ let us take $J^{(k)}_{i}=G_{k}\cap\ac{i+1,\ldots,n}$, $\ell^{(k)}=m_k$, $\tau=1/14$, $r_{0}=2$ and 
$S_i^{(k)}= \sum_{b\in J^{(k)}}U_{ib}=\sum_{b\in J^{(k)}}E_{ib}$. Define accordingly, the sets   $I_{r}^{(k)}$ and 
\begin{equation}\label{eq:definition_Omega_k}
\Omega_{k}=\ac{\sum_{i\in I^{(k)}_{0}} (S_{i}^{(k)})^2 \al nm_{k} L}\bigcap \bigcap_{r\geq 1} \ac{ |I^{(k)}_{r}|\leq 2n 2^{-\tau m_{k} L r2^{r_{0}+r}}}\ ,
\end{equation}
and $y^{(k)}_{r}=2^{r_{0}+r} \ell^{(k)} L$.
Then, according to Lemma~\ref{lem:peel}, the event 
$\Omega=\cap_{k=1}^K\Omega_{k}$ 
holds with probability at least $1-c/n^2$.
Let us now prove (\ref{sum:ordered:sbm}). 
We consider apart the case $j\neq k$ and $j=k$. In the remainder of the proof, $k$ is fixed and to alleviate the notation, we simply write $I_r$ and $S_i$ for $I_{r}^{(k)}$ and $S_i^{(k)}$. 
\medskip

\noindent{\bf Case $j\neq k$.}
Let $Q$ be a subset of $G_{j}$ with cardinality $q$, and set $T_{i,Q}=\sum_{b\in Q} U_{ib}$ for $i\in \{1,\ldots, n\}$. Then, we have 
\[
 \sum_{b\in Q}z^{(j,k)}_{b}= \frac{1}{m_k} \sum_{i=1}^n  S_iT_{i,Q}\ .
\]
Since all the entries of $U$ are independent and $Q\cap J_i=\emptyset$, the $T_{i,Q}$'s are independent from the $S_i$'s. 
Let us first upper bound the sum
$\sum_{i\in I_{0}} S_{i}T_{i,Q}=\sum_{i=1}^n (S_{i}\1_{S_{i}^2\leq y^2_{1}})T_{i,Q}$
on $\Omega$. Working conditionally on the $(S_{i}:i\in I_{0})$, we have 
with Bernstein inequality
$$\P\cro{\sum_{i\in I_{0}} S_{i}T_{i,Q} \ag \sqrt{qL\sum_{i\in I_{0}}S_{i}^2 t}+m_{k}Lt \bigg|S_{i}:i\in I_{0}}\leq e^{-t}\ ,
$$
Hence, since $\sum_{i\in I_{0}}S_{i}^2\al nm_{k}L$ on $\Omega$, we have
\begin{equation}\label{eq:upper1:I0}
\P\cro{\Omega\cap \sum_{i\in I_{0}} S_{i}T_{i,Q} \ag  L\sqrt{qnm_{k}t}+m_{k}Lt}\leq e^{-t} \ ,
 \end{equation}
for any $t>0$. Let us now upper-bound $\sum_{r\geq 1}\sum_{i\in I_{r}} S_{i}T_{i,Q}$ on $\Omega$. 
Since such $S_{i}$ are positive, we have for $\lambda>0$
\begin{align*}
\E\cro{\exp\pa{\lambda\sum_{r\geq 1}\sum_{i\in I_{r}}S_{i}T_{i,Q}}\1_{\Omega}} &\leq \E\cro{\exp\pa{Lq\sum_{r\geq 1}\sum_{i\in I_{r}}(e^{\lambda S_{i}}-1-\lambda S_{i})}\1_{\Omega}}\\
&\leq \E\cro{\exp\pa{Lq\sum_{r\geq 1} |I_{r}|(e^{2\lambda m_{k}L 2^{r_{0}+r}}-1-2\lambda m_{k}L 2^{r_{0}+r})}\1_{\Omega}}\\
&\leq \exp\pa{2nLq\sum_{r\geq 1} 2^{-\tau m_{k}L r2^{r_{0}+r}}e^{2\lambda m_{k}L 2^{r_{0}+r}
}}
\end{align*}
For $\lambda=\tau \log(2)/4$, we have
\begin{align*}
 \E\cro{\exp\pa{{\tau \log(2)\over 4}\sum_{r\geq 1}\sum_{i\in I_{r}}S_{i}T_{i,Q}}\1_{\Omega}} & \leq  \exp\pa{2nLq\sum_{r\geq 1} 2^{-\tau m_{k}L r2^{r_{0}+r-1}}}
 \\
 &\leq e^{2nLq e^{-cm_{k}L}} \leq e^{nLq} \ ,
\end{align*}
since $L\gtrsim m_k$. This gives 
$$\P\cro{\Omega \cap \sum_{r\geq 1}\sum_{i\in I_{r}}S_{i}T_{i,Q} >t}\leq e^{-{\tau \log(2) t\over 4}+nLq}.$$
Together with \eqref{eq:upper1:I0}, we obtain that 
$$\P\cro{\Omega \cap \sum_{i=1}^nS_{i}T_{i,Q} \ag nLq+m_{k}Lt}\leq 2e^{-t}$$
and so 
$$\P\cro{\Omega \cap \max_{Q\subset G_{j},\, |Q|= q}\sum_{b\in Q}z^{(j,k)}_b\ag {nLq\over m_{k}}+Lt}\leq  2C_{m_{j}}^qe^{-t}.$$

\bigskip
\noindent{\bf Case $j=k$.} We start from $\sum_{b\in Q}z^{(j,k)}_b= - m_k^{-1}\sum_{i=1}^n S_iT_{i,Q}$. Unfortunately, the sums  $S_i$ and $T_{i,Q}$  are no longer independent as $Q\subset G_k$ and $T_{i,Q}$ is therefore a subsum of $S_i$. 
We define $k(i)$ as the index in $\ac{1,\ldots,K}$ such that $i\in G_{k(i)}$ and we
set $J_{i}=G_{k}\cap\ac{i+1,\ldots,n}$, $L_{i}=P_{k(i)k}$ and  $N_{i}=\sum_{a\in J_i}X_{ia}=S_{i}+|J_{i}|L_{i}$.
We observe that conditionally on $N_{i}$, the sum $H_{i}=\sum_{b\in Q\cap \ac{i+1,\ldots,n}}X_{ib}$ follows a Hypergeometric distribution with parameter $(q_{i},N_{i},|J_{i}|)$, where $q_{i}=|Q\cap\ac{i+1,\ldots,n}|\leq q$. Let $H'_{i}$ be a random variable with binomial $(q_{i},N_{i}/|J_{i}|)$ distribution conditionally on $N_{i}$.
Since $x\to e^{-\lambda S_{i}x}$ is continuous and convex, according to Theorem 4 of \cite{Hoeffding}, we have conditionally on $S_{i}$
$$\E[e^{-\lambda S_{i}H_{i}}|S_{i}]\leq \E[e^{-\lambda S_{i}H_{i}'}|S_{i}]\ .$$
Hence, conditionally on $S_{i}$, we apply Chernoff bound to $\sum_{i\in I_{0}}S_{i}H_i$ together with the above control of the Laplace transform. This allows us to get a Bernstein like inequality. Hence, with probability higher than $1-e^{-t}$ we have 
\begin{align*}
-\sum_{i\in I_{0}}S_{i}T_{i,Q}&=-\sum_{i\in I_{0}}S_{i}\pa{H_i-q_{i}L_{i}}\\
&\leq -\sum_{i\in I_{0}}q_{i}S_{i}\pa{{N_{i}\over |J_{i}|}-L_{i}}+ c \sqrt{\sum_{i\in I_{0}}q_iS_{i}^2 {N_{i}\over |J_{i}|} t}+c m_{k}Lt\\ 
&\leq - \sum_{i\in I_{0}}{q_i S_{i}^2\over |J_{i}|}+ c \sqrt{\max_{i\in I_0}N_i \sum_{i\in I_{0}}{q_i S_{i}^2\over |J_{i}|} t}+c m_{k}Lt \\
&\al  m_{k}Lt,
\end{align*}
where we used in the last line that  $N_{i}=S_{i}+|J_{i}|L_{i} \leq (1+2^{r_{0}})m_{k}L$ for $i\in I_{0}$.

When $S_{i}>0$, we have 
$$-S_{i}T_{i,Q} \leq qLS_{i}\ ,$$
since each entry $E_{ab}$ is larger or equal to $-L$. 
For 
 $i\in I_{r}$  with $r\geq 1$, we have $y_{r}<S_{i}\leq 2 y_{r}$ and hence, on $\Omega$
\begin{align}
-\sum_{i}\sum_{r\geq 1} S_{i}T_{i,Q}\1_{i\in I_r} & \leq 2qL \sum_{r\geq 1} |I_{r}| y_{r}\nonumber \\
&\leq 4nqL \sum_{r\geq 1} 2^{r_{0}+r} m_{k}L 2^{-\tau m_{k}L r 2^{r_{0}+r}}\nonumber \\
&\leq 4nqL 2^{-cm_{k}L}\ , \nonumber
\end{align}
since $m_kL\gtrsim 1$. 
Hence arguing as before we obtain that for $t\geq 0$
$$\P\cro{\Omega\cap \max_{Q\subset G_{k}:|Q|= q}\sum_{b\in Q}z^{(j,k)}_b\ag {nLq\over m_{k}}+Lt}\leq  C_{m_{j}}^qe^{-t}.$$
The proof of Lemma~\ref{lem:quadratic:sbm} is complete.

\subsubsection{Case $(U^T)^2$}
The case $(U^T)^2$ is somewhat more messy, due to interlaced rows/columns dependences. Recall, that for  two indices $j$ and $k$, $(-1)^{j=k}$ equals $-1$ when $i=k$ and $1$ when $i\neq k$. For $1\leq j\leq k$ and $b\in G_j$, define 
\[
z^{(j,k)}_b = {\pa{-1}^{j=k}\over m_{k}} \sum_{i=1}^n\sum_{a\in G_{k}} U_{ia}U_{bi}\ , 
\]
and let $z^{(j,k)}_{(1)}\geq z^{(j,k)}_{(2)}\geq \ldots$ be the random variables $z^{(j,k)}_b$ ranked in decreasing order. Compared to $U^TU$, the difficulty is that the same random variables $U_{ia}$ can occur several times in the definition of $z^{(j,k)}_b$.

\begin{lem}\label{lem:square:sbm}
There exists an event $\underline{\Omega}$  of probability at least $1-c/n^2$, such that for all $1\leq j,k\leq K$ and  for any integer $q \in [1,m_{k}]$ and $t\geq 0$, we  have 
 \begin{equation}\label{sum:ordered-square:sbm}
 \P\cro{\underline{\Omega}\cap\ac{\sum_{u=1}^q z_{(u)}^{(j,k)}\ag L\pa{{qn\over m_k}+t}}}\leq c''C_{m_{k}}^qe^{-c't}.
 \end{equation}
\end{lem}
As for the previous case, we easily conclude from Lemma \ref{lem:square:sbm} by applying a union bound together with Jensen inequality.

\medskip 

\noindent{\bf Proof of Lemma~\ref{lem:square:sbm}.} 
The event  $\underline{\Omega}$ is defined as the intersection $\underline{\Omega}=\Omega\cap \cap_{u=1}^7\cap_{k=1}^K \underline{\Omega}^{(u)}_k$, where the event $\Omega=\cap_{k}\Omega_k$ is introduced in \eqref{eq:definition_Omega_k} the previous proof and the events  $\underline{\Omega}_k^{(u)}$ are defined along the proof. Let us split the sum $\pa{-1}^{j=k} \sum_{i=1}^n\sum_{a\in G_{k}, b\in Q} U_{ia}U_{bi}$
into two parts depending whether $i\in G_k$ or not.

\medskip 

\noindent {\bf Case $i\notin G_{k}$.} 
For $Q\subset G_{j}$, consider the sum 
\begin{equation}\label{eq:_sum_SiS'i}
\pa{-1}^{j=k} \sum_{i\notin G_{k}}\sum_{a\in G_{k}, b\in Q} U_{ia}U_{bi}=\pa{-1}^{j=k} \sum_{i\notin G_{k}} S_{i}T_{i,Q}'\ , 
\end{equation}
with $S_{i}=\sum_{a\in G_{k}} U_{ia}$ as in the proof of Lemma \ref{lem:quadratic:sbm}  (we dropped the exponent $(k)$ to alleviate the notation) and $T_{i,Q}'=\sum_{b\in Q} U_{bi}$. In the collection of $(S_i)$'s and $(T_{i,Q}')$'s, all the random variables
are independent since the sums respectively run on the sets $G_{k}^c\times G_{k}$ and $Q\times G_{k}^c$ that do not intersect. As a consequence, \eqref{eq:_sum_SiS'i} is handled exactly as the case $j\neq k$ in Lemma~\ref{lem:quadratic:sbm}. We conclude that for all $1\leq j,k\leq K$ and all $Q\subset G_j$ of size $q$, we have 
\begin{equation}
 \label{eq:lemme_horrible_conclusion1}
\P\cro{ \Omega \cap  \left\{\pa{-1}^{j=k} \sum_{i\notin G_{k}}\sum_{a\in G_{k}, b\in Q} U_{ia}U_{bi}\ag {nLq\over m_{k}}+Lt\right\}}\al  e^{-t}\ .
\end{equation}

\medskip 
\noindent{\bf Case $i\in G_k$ and $j\neq k$.} For $Q\subset G_j$, we consider the sum 
\begin{equation}\label{eq:_sum_SiS'i_2}
\pa{-1}^{j=k} \sum_{i\in G_{k}}\sum_{a\in G_{k}, b\in Q} U_{ia}U_{bi}=\pa{-1}^{j=k} \sum_{i\notin G_{k}} S_{i}T_{i,Q}'\ , 
\end{equation}
with $S_{i}=\sum_{a\in G_{k}} U_{ia}$ and $T_{i,Q}'=\sum_{b\in Q} U_{bi}$. As above, the indices run in  $G_{k}\times G_{k}$ and $Q\times G_{k}$ which do not intersect.  Again, \eqref{eq:_sum_SiS'i_2} is handled exactly as the case $j\neq k$ in Lemma~\ref{lem:quadratic:sbm}.  We conclude that for all $1\leq j,k\leq K$ with $j\neq k$ and all $Q\subset G_j$ of size $q$, we have 
\begin{equation}
 \label{eq:lemme_horrible_conclusion2}
\P\cro{ \Omega \cap \left\{  \sum_{i\notin G_{k}}\sum_{a\in G_{k}, b\in Q} U_{ia}U_{bi}\ag {nLq\over m_{k}}+Lt\right\}}\al  e^{-t}\ .
\end{equation}

\medskip 

\noindent {\bf Case $i\in G_{k}$ and $j=k$.}  
It remains to upper-bound the sum
\begin{equation}\label{eq:equation_horrible}
-\sum_{i\in G_{k}} \sum_{a\in G_{k}} U_{ia}\sum_{b\in Q} U_{bi}\ ,
\end{equation}
for $Q\subset G_{k}$.
It is the main hurdle of the proof. Indeed, we multiply row sums  $\sum_{a\in G_{k}} U_{ia}$ of the matrix $U$ (restricted to $G_k\times G_k$) to columns sums of the matrix $U$. As in the previous proof, we consider separately small and large row sum of $U$.
Define the set of indices corresponding to small $S_i$'s
\begin{equation}\label{eq:definition_I0_alternative_horrible}
\underline{I}_0= \{i \in G_k:  S_i\leq  5 m_k L\}\ , 
\end{equation}
where we recall that  $S_{i}=\sum_{a\in G_{k}} U_{ia}$. 

Define 
$\ell_{b}=|G_{k}\cap \ac{b+1,\ldots,n}|\leq m_{k}$ and $N_{b}=S_{b}+\ell_{b} P_{kk}$. Given the collection $(S_i)$, $i\in G_k$, the binary random variable $U_{bi}+P_{kk}$ is distributed as a sampling of size 1 in an urn of size $\ell_b$ containing $N_b$ ones. Hence, we split the sum \eqref{eq:equation_horrible} into three  pieces to center the random variables $U_{bi}$.  
\begin{align*}
\lefteqn{-\sum_{i\in G_k}  \sum_{a\in G_k, a > i} \sum_{b\in Q, b< i}  U_{ia}U_{bi} }\\
&=  - \sum_{b\in Q}\sum_{i\in G_k, \ i > b} S_i  U_{bi} \\
&  =  - \sum_{i \in G_k} S_i \sum_{b\in Q,\ b <i } \left( \frac{N_b}{\ell_b}-P_{kk}\right)  - \sum_{i \in G_k} S_i \sum_{b\in Q,\ b <i }\left(U_{bi}+P_{kk}-  \frac{N_b}{\ell_b} \right) \\
& \leq    - \sum_{b\in Q} \left( \frac{N_b}{\ell_b}-P_{kk}\right) \sum_{i \in G_k, \ i >b} S_i     - \sum_{i \in \underline{I}_0} S_i \sum_{b\in Q,\ b <i }\left(U_{bi}+P_{kk}-  \frac{N_b}{\ell_b} \right)  +  \sum_{i \in \underline{I}_0^c} S_i \sum_{b\in Q,\ b <i } \frac{N_b}{\ell_b}  \\ 
& = E_1+ E_2 + E_3\ ,
\end{align*}
where we used that $S_i> 0 $ for $i\in \underline{I}_0^c$. 

We shall prove that the three following bounds  hold, for any $t\geq 1$, 
\begin{eqnarray}	
\label{eq:control_E1}
1_{\underline{\Omega}}|E_{1}|\lesssim qnL\ ,  \\
\label{eq:E2_conclusion}
\P\cro{\underline{\Omega} \cap \{E_{2}\ag qnL + m_kLt\}} &\leq& 2 e^{-t}\ ,\\
\label{eq:E3}
 \P\cro{\underline{\Omega} \cap \{E_{3} \ag qnL + t\}}& \leq& e^{-t}\ ,
\end{eqnarray}
and that $\P[\underline{\Omega}]\geq 1-c/n^2$.
Gathering these three bounds, we obtain for all $1\leq k\leq K$ and all $Q\subset G_k$, we have 
\begin{equation}
 \label{eq:lemme_horrible_conclusion3}
\P\cro{  - \1_{\Omega}\sum_{i\notin G_{k}}\sum_{a\in G_{k}, b\in Q} U_{ia}U_{bi}\ag {nLq\over m_{k}}+Lt}\al  e^{-t}\ .
\end{equation}
Together with \eqref{eq:lemme_horrible_conclusion1} and \eqref{eq:lemme_horrible_conclusion2}, this concludes the proof.  It remains to show~(\ref{eq:control_E1}--\ref{eq:E3}).

 \medskip

\noindent {\bf Control of $E_1$}. The random variable $\sum_{i>b,\ i\in G_k}S_i$ is distributed as a Binomial random variable  with parameters $P_{kk}$ and $n_{b}\leq \ell_b^2/2$. 
By Bernstein inequality together with a union bound, we derive that,  on an event $\underline{\Omega}_k^{(1)}$ with probability higher than $1-1/n^3$, we have 
\begin{eqnarray} \label{eq:Bernstein1}
| \sum_{i>b}S_i| &\lesssim& \ell_b \sqrt{ L\log(n) }+ \log(n)\\
|N_b| &\lesssim&\ell_{b} L + \sqrt{ \ell_b L\log(n) }+ \log(n)\ \al \ell_{b} L+\log(n),  \label{eq:Bernstein2}
\end{eqnarray}
 uniformly on $b\in G_k$.  Consider any $b$ such that $\ell_b^2L\geq {\log(n)}$. Then, 
 \begin{eqnarray}\label{eq:upperE1-1}
 \big|\left(\frac{N_b}{\ell_b} -P_{kk}\right)\sum_{i \in G_k, \ i >b} S_i \big| &\lesssim&   L^{3/2} {\ell_b} \sqrt{\log(n)} + \log^{3/2}(n)\sqrt{L} \lesssim nL\ , 
 \end{eqnarray}
since $m_kL\geq 1$ and  $L\leq 1/\log(n)$. 

Now assume that $\ell_b^2 L\leq \log(n)$ and $L\geq \frac{\log(n)}{n}$. By definition, we have $N_b\leq \ell_b$ almost surely. Together with \eqref{eq:Bernstein1}, this leads us to
\begin{eqnarray}\label{eq:upperE1-2-0}
 \bigg|\left(\frac{N_b}{\ell_b}-P_{kk}\right) \sum_{i \in G_k, \ i >b} S_i \bigg| &\lesssim&   \log(n) \leq nL\ . 
\end{eqnarray}
Next, we consider the case where $\ell_b^2 L\leq \log(n)$ and   $L\leq \frac{\log(n)}{n}$. If $\ell_b> n^{1/4}$, then \eqref{eq:Bernstein1} leads us to 
 \begin{eqnarray}\label{eq:upperE1-2}
 \bigg|\left(\frac{N_b}{\ell_b}-P_{kk}\right)\sum_{i \in G_k, \ i >b} S_i \bigg| &\lesssim& L\log(n) + \frac{\log^2 (n)}{n^{1/4}}   \lesssim 1 \lesssim nL 
 \end{eqnarray}
 Let $b_0$ such that $\ell_{b_{0}}= \lfloor n^{1/4}\rfloor$ (if it exists). Then, $ \sum_{i \in G_k, \ i >b_0} N_i$ follows a  binomial distribution with parameters $P_{kk}\leq \log(n)/n$ and  $r\leq \sqrt{n}$. On an event $\underline{\Omega}_k^{(2)}$ with  probability higher than $1- 1/n^3$, it is no higher than $5$. Under this event, we have for any $b\geq b_0$, 
\begin{equation}\label{eq:upperE1-3}
 \bigg|\left(\frac{N_b}{\ell_b}-P_{kk}\right) \sum_{i \in G_k, \ i >b} S_i \bigg| \lesssim    1 \leq nL\ . 
 \end{equation}
 Gathering (\ref{eq:upperE1-1}--\ref{eq:upperE1-3}) and summing over all $b\in Q$,  we have proved \eqref{eq:control_E1}.

\bigskip

\noindent {\bf Control of $E_2$}. We work conditionally to $S_i$. In such a case, for a fixed $b$, the random variables  $((U_{bi}+P_{kk})_{i>b})$ are distributed as a sampling without replacement of $N_b$ ones in an urn of size $\ell_b$. Then,  according to Theorem 4 of \cite{Hoeffding}, the Laplace transform of
$-\sum_{b\in Q} \sum_{i> b}  S_i  (U_{bi}+P_{kk}-N_b/\ell_b)$ conditional to the $(S_{i})_{i}$ is upper bounded by that of $-\sum_{b\in Q}\sum_{i>b}S_i  (\tilde{E}_{ib}-N_b/\ell_b)$ where the $\tilde{E}_{ib}$ are independent and follow a Bernoulli distribution with parameters $N_b/\ell_b$. 

Hence, we can apply Bernstein's inequality conditionally to $S_i$ to obtain that with probability at least $1-e^{-t}$ 
\begin{eqnarray}\label{eq:E2_0}
 E_2 &\lesssim &   \sqrt{\sum_{i\in \underline{I}_0} S_i^2 \bigg( \sum_{b\in Q,\ b<i }\frac{N_b}{\ell_b}\bigg) t } + m_kL t\,, 
\end{eqnarray}
since $\sup_{i\in \underline{I}_0}|S_i|\leq 5 m_k L$.
We define $s^*=\log^2(n)$, $Q_{-}=\ac{b\in Q: \ell_{b}\leq s^*}$ and $Q_{+}=\ac{b\in Q: \ell_{b}> s^*}$.
Then, we split the sum into two parts
\[
\sum_{i\in \underline{I}_0} S_i^2  \sum_{b\in Q,\ b<i }\frac{N_b}{\ell_b} = \sum_{i\in \underline{I}_0} S_i^2  \sum_{b\in Q_{+},\ b<i }\frac{N_b}{\ell_b}+ 
 \sum_{i\in \underline{I}_0} S_i^2 \sum_{b\in Q_{-},\ b<i }\frac{N_b}{\ell_b}
\]

\noindent 
{\bf Sum over $Q_{+}$}. By Bernstein inequality, we have, simultaneously for all $i\in G_k$, 
\[
 |S_i|\1_{i\in \underline{I}_0}\lesssim \sqrt{m_kL\log(n)} + \left( \log(n) \wedge (m_kL)\right) \lesssim \sqrt{m_kL \log(n)}
\]
on an event $\underline{\Omega}^{(3)}_k$ with probability higher than $1-1/n^3$. Since the random variables $S_i^2\1_{i\in \underline{I}_0}$ are independent and their variance is less than $6\ell_{i}^2L^2+ \ell_{i}L$, we derive from Bernstein inequality that, on an event $\underline{\Omega}^{(3)}_k\cap\underline{\Omega}^{(4)}_k$ with probability higher than $1-2/n^3$,  for all  $b\in Q_{+}$, we have 
\begin{eqnarray}
 \sum_{i>b} S_i^2\1_{i\in \underline{I}_0} &\lesssim &\ell_b^2 L + \sqrt{[\ell_b^3L^2 + \ell_b^2 L]\log(n)} + m_kL \log^2(n) \nonumber \\ 
 &\lesssim& \ell_bm_k  L + m_kL \log^2(n)\lesssim  \ell_bm_kL\ ,\label{eq:E2_1} 
 \end{eqnarray}
since $\log^2(n)\leq \ell_b\leq m_k$ and $m_kL \geq 1$. Under this event, we obtain 
\begin{eqnarray*}
\sum_{i\in \underline{I}_0} S_i^2  \sum_{b\in Q_{+},\ b<i }\frac{N_b}{\ell_b} &\al   & m_kL \sum_{b\in Q_{+}}N_b
\end{eqnarray*}
Then, $\sum_{b\in Q_{+}}N_b$ is stochastically dominated by a Binomial distribution with parameters $L$ and $qm_k$, hence we derive from Bernstein inequality that
$$ \sum_{b\in Q_{+}}N_b \al qm_{k}L+t\ ,$$
with probability higher than $1-e^{-t}$. Gathering \eqref{eq:E2_0}, \eqref{eq:E2_1}, and the last bound, we obtain 
\begin{equation}\label{eq:E2_3}
\P\pa{\underline{\Omega} \cap \{E_{2}^+\ag m_{k}L(\sqrt{qt}+t)\}} \leq 2 e^{-t},
\end{equation}
where $E_{2}^+$ corresponds to the sum $E_{2}$ restricted to the indices $b\in Q_{+}$.

\bigskip 

\noindent 
{\bf Sum over $Q_{-}$}. We consider two subcases. First assume that $L\geq n^{-1/4}$. Since $\sum_{i: \ell_{i}\leq s^*}S_i^2\leq \sum_{i=1}^{\log^2(n)-1} i^2  \leq \log^6(n)$ and, since $N_b\leq \ell_{b}$, we derive that 
\begin{equation}\label{eq:E2_4}
 \sum_{i\in \underline{I}_0} S_i^2  \sum_{b\in Q_{-},\ b<i }\frac{N_b}{\ell_b} \leq    q\log^6(n)\lesssim qL^2nm_k\ . 
\end{equation}
Next, we assume that $L\leq n^{-1/4}$. Let $b_{*}=\min\ac{b\in G_{k}: \ell_{b}\leq s^*}$. Since, by definition of $s^*$ and $b_{*}$, we have
$$\big| \ac{(i,a):i,a\in G_{k},\ i\geq b_{*}+1,\ a\geq i+1}\big| \leq (s^*)^2/2= \log^4(n)/2,$$
the sums $(S_i)_{i>b_{*}}$ involves less than $\log^{4}(n)/2$ independent Bernoulli random variables with parameters less that $L$.
Hence, on an event $\underline{\Omega}^{(5)}_k$ with probability larger than $1- 1/n^{3}$, at most $10$ of them are equal to one and
\begin{eqnarray}\label{eq:E2_5}
 \sum_{i\in \underline{I}_0} S_i^2  \sum_{b\in Q_{-},\ b<i }\frac{N_b}{\ell_b}\lesssim q(L^2 \log^6(n)+1) \lesssim qL^2nm_k\ , 
\end{eqnarray}
since $m_kL\geq 1$.
Gathering (\ref{eq:E2_0}--\ref{eq:E2_5}), we have proved \eqref{eq:E2_conclusion}.

\medskip 

\noindent {\bf Control of $E_3$}.
If $L\geq \log(n)/m_k$, then Bernstein inequality enforces that $\underline{I}_0^c= \emptyset$ and therefore $E_3=0$ with probability higher than $1-1/n^3$. Let us call $\underline{\Omega}^{(6)}_k$ the corresponding event.  Hence, we  assume henceforth that $L\leq \log(n)/m_k$. We claim that, on an event  $\underline{\Omega}^{(7)}_k$ with probability larger than $1-1/n^3$, we have
 \begin{equation} \label{eq:upper_tail_Si} 
  \sum_{i>b} S_i\1_{i\in \underline{I}_0^c}\lesssim \ell_b+\log^3(n)\ ,
 \end{equation}
  uniformly over all $b\in G_{k}$.  The proof of this claim is a slight variation on the proof of Lemma~\ref{lem:peel}. We provide it here for the sake of completeness.  With probability higher than $1-1/n^3$, we have 
$\max_{i}S_i\lesssim \log(n)$.   Write $\log_2$ for the binary logarithm. Fix any $b\in G_k$ and decompose 
 \begin{equation*}
  \sum_{i>b} S_i\1_{i\in \underline{I}_0^c}\leq \sum_{r=1}^{\lfloor \log_2(c''\log(n)/m_kL)\rfloor} \sum_{i>b} 5 2^r L m_k \1_{S_i\in [5\cdot 2^{r-1} L m_k, 5\cdot2^{r} L m_k]}
 \end{equation*}
The random variables  $\sum_{i>b}  \1_{S_i\in [ 5\cdot 2^{r-1} L m_k,  5\cdot 2^{r} L m_k]}$ are stochastically dominated by  binomial distributions with parameters $\ell_b$ and $p_{r}$, where $p_r\leq e^{-c' 2^r Lm_k}$ is the probability that a Binomial distribution with parameters $(m_k,L)$ is larger than $ 5\cdot 2^{r-1} L m_k$. Applying Bernstein's inequality together  with a union bound, we conclude that, simultaneously for all $b$ and all $r$, 
\[
\sum_{i>b}  \1_{S_i\in [ 5\cdot 2^{r-1} L m_k,  5\cdot2^{r} L m_k]}\lesssim \ell_b e^{-c' 2^{r-1} Lm_k}+\log(n)\ ,
\]
with probability higher than $1-1/n^3$.
This leads us to 
 \begin{equation*}
  \sum_{i>b} S_i\1_{i\in \underline{I}_0^c}\leq \sum_{r=1}^{\lfloor \log_2(c''\log(n)/m_kL)\rfloor}   c 2^r L m_k \left[ \ell_b e^{-c' 2^r Lm_k}+\log(n)\right]\lesssim \ell_b+ \log^3(n)\ ,
 \end{equation*}
 since $Lm_k\geq 1$. We have proved the claim \eqref{eq:upper_tail_Si}.

 \bigskip 
 
 As for $E_{2}$, we decompose $Q=Q_{+}\cup Q_{-}$, with $s^*$ now set to $s^*=\log^3(n)$.

 \noindent 
{\bf Sum over $Q_{+}$}. 
Let us work on the event $\underline{\Omega}^{(7)}_k$ so that  \eqref{eq:upper_tail_Si} holds. Hence, as for $E_{2}$, according to Bernstein inequality, we have with probability larger than $1-e^{-t}$
 \begin{equation}\label{eq:E3_1}
  \sum_{b\in Q_{+}}   \frac{N_b}{\ell_b}  \sum_{i \in \underline{I}_0^c,\ i >b} S_i \lesssim       \sum_{b\in Q_{+}}  N_b \al qm_{k}L+t.
  \end{equation}

\noindent 
{\bf Sum over $Q_{-}$}. As for $E_2$, we consider two subcases depending whether $L\leq n^{-1/4}$ or $L> n^{-1/4}$. If $L\leq n^{-1/4}$, we have, as argued in  $E_2$, that $\sum_{i>b_{*}} N_{i}$ is less than $10$ on an event $\underline{\Omega}^{(5)}_k$ with probability larger than $1- 1/n^{2}$. Under this event, we have 
\begin{equation}\label{eq:E3_2}
 \sum_{b\in Q_{-}}   \frac{N_b}{\ell_b}  \sum_{i \in \underline{I}_0^c,\ i >b} S_i \lesssim   q  \lesssim q m_k L \,.
 \end{equation}
 If $L>n^{-1/4}$, we  use the crude bound 
\begin{equation}\label{eq:E3_3}
  \sum_{b\in Q_{-}}   \frac{N_b}{\ell_b}  \sum_{i \in \underline{I}_0^c,\ i >b} S_i \leq    q (s^*)^2 \lesssim qL n\ . 
 \end{equation}
 Putting (\ref{eq:E3_1}--\ref{eq:E3_3}) together, we conclude that on the event $\underline{\Omega}$ which  has a probability larger than $1-7/n^2$ we have proved \eqref{eq:E3}.

\subsection{Proof of Lemma~\ref{lem:peel}}
With no loss of generality we only consider the case $I=\ac{1,\ldots,n}$. 
Let us first handle the upper-bound on the sum $\sum_{i=1}^n S_{i}^2\1_{S_{i}^2\leq y^2_{1}}$, with $y_{1}=2^{r_{0}+1}\ell L$. We observe that the variance of $S_{i}^2$ is upper bounded by $3\ell L$, hence,
according to Bernstein inequality,
$$\sum_{i=1}^nS_{i}^2\1_{S_{i}^2\leq y^2_{1}} \al \ell L(n+ \sqrt{\log(n)})+y_{1}^2\log(n),$$
with probability at least $1-1/n^3$.
We observe that $$y_{1}^2\log(n)\al (\ell L)^2\log(n)\al n\ell L$$ according to the assumptions on $L$, so
with probability at least $1-1/n^3$ 
$$\sum_{i=1}^nS_{i}^2\1_{S_{i}^2\leq y_{1}^2}\al n\ell L.$$

We now turn to the first part of the lemma. We shall first work around a bound of $\P(y_{r}< S_{i} \leq y_{r+1})$.
We set $h(x)=(1+x)\log(1+x)-x\geq x\log(x/4)$ for $x\geq 0$ and $y_{r}=2^{r_{0}+r}\ell L$. Denoting $\sigma^2_i$ the variance of $S_i$, we deduce from 
Bennett's inequality that, for $r\geq 1$,
\begin{eqnarray*}
\P(y_{r}< S_{i} \leq y_{r+1})\leq \P( S_{i}>y_r)&\leq& e^{- \sigma^2_ih(y_r/\sigma^2_i)}\leq \exp\left[-y_r \log\left(\frac{y_r}{4\sigma^2_i}\right)
\right]\\
&\leq & 
2^{-\ell L (r+r_{0}-2)2^{r_{0}+r}}=: p_{r}\  , 
\end{eqnarray*}
since $\sigma_i^2\leq \ell L$. 

Next, we use again Bennett inequality to ensure that $|I_r|$ is small. 
For $1/14\leq \tau\leq 1/2$, and $r\geq 1$, Bennett's inequality  again ensures that
$$\P[|I_{r}|\geq 2np_{r}^{\tau}]\leq \exp(-np_{r} h(p_{r}^{-(1-\tau)}))\leq \exp(-0.5 np_{r}^\tau \log(p_{r}^{-(1-\tau)}))\ ,$$
since $h(x)\geq x\log(x)/2$ for $x\geq e^{2}$ and since  $\tau\leq 1/2$ and $p_r\leq 2^{-16}$ is small enough. 
Since $r_{0}\geq 2$, we have $2np_{r}^{\tau}< 1$ for $r\geq r^*= \log_{2}(\tau^{-1}\log(2n))$. 
Let $\Omega'$ be the event $$\Omega'=\bigcap_{r\geq 1} \ac{|I_{r}|\leq 2np_{r}^{\tau}}=\bigcap_{r=1}^{\lfloor r^*\rfloor } \ac{|I_{r}|\leq 2np_{r}^{\tau}}\bigcap \ac{\big|\{i: S_i>y_{r^*}\}\big|=0}  .$$
So we have
\[
\P(\Omega'^c)\leq    \sum_{r=1}^{\lfloor r^*\rfloor }  \exp(-0.5 np_{r}^\tau \log(p_{r}^{-(1-\tau)}))+   \exp(-0.5 np_{r^*}^\tau \log(p_{r^*}^{-(1-\tau)}))  \ . 
\]
 For $x\in (0,1)$, the function $\phi_{\tau}:x\mapsto x^{\tau}\log(x)$ is decreasing when $\tau\log(x)<-1$ and then increasing. Since $(2n)^{-1/\tau}\leq p_{r} \leq 2^{-16}$, $\phi_{\tau}(p_r)$ achieves its maximum at $p_r= (2n)^{-1/\tau}$ for $\tau \in [1/10,1/2]$. For $\tau\in [1/14,1/10)$, its maximum its achieved either at $p_r= 2^{-16}$ or $p_r= (2n)^{-1/\tau}$ and it can be checked to be reached at $p_r= (2n)^{-1/\tau}$ as long as $(2n)^{-1/\tau}\leq 2^{-16}$. Note that for smaller $nT$, the event $\Omega'$ reduces to $\ac{\big|\{i: S_i>y_{r^*}\}\big|=0}$.  
 For $n\geq 2$, we conclude that 
\[
 \P(\Omega'^c)\leq (\lfloor r^*\rfloor +1) \exp\pa{- {1-\tau\over 4\tau} \log(2n)} \leq  {\log_{2}(\tau^{-1}\log(2n))+1 \over (2n)^{(1-\tau)/4\tau}}\, .
\]

\section*{Acknowledgements}
We thank the editors and anonymous reviewers for their helpful suggestions. 
 Christophe Giraud is partially supported by the LabEx LMH, ANR-11-LABX-0056-LMH and the CNRS PICS grant HighClust.

\bibliographystyle{alpha}
\bibliography{biblio}

\newcommand{\etalchar}[1]{$^{#1}$}
\begin{thebibliography}{BGRV16}

\bibitem[{Abb}17]{AbbeReview2017}
E.~{Abbe}.
\newblock {Community detection and stochastic block models: recent
  developments}.
\newblock {\em ArXiv e-prints}, March 2017.

\bibitem[ABC{\etalchar{+}}15]{awasthi2015relax}
Pranjal Awasthi, Afonso~S Bandeira, Moses Charikar, Ravishankar Krishnaswamy,
  Soledad Villar, and Rachel Ward.
\newblock {Relax, no need to round: Integrality of clustering formulations}.
\newblock In {\em {Proceedings of the 2015 Conference on Innovations in
  Theoretical Computer Science}}, pages 191--200. ACM, 2015.

\bibitem[ACKS15]{awasthi2015hardness}
Pranjal Awasthi, Moses Charikar, Ravishankar Krishnaswamy, and Ali~Kemal Sinop.
\newblock {The Hardness of Approximation of Euclidean k-means}.
\newblock {\em arXiv preprint arXiv:1502.03316}, 2015.

\bibitem[AL14]{amini}
A.~A. {Amini} and E.~{Levina}.
\newblock {On semidefinite relaxations for the block model}.
\newblock {\em ArXiv e-prints}, June 2014.

\bibitem[AM05]{AchlioptasMcSherry2005}
Dimitris Achlioptas and Frank McSherry.
\newblock {On Spectral Learning of Mixtures of Distributions}.
\newblock In Peter Auer and Ron Meir, editors, {\em {Learning Theory}}, pages
  458--469, Berlin, Heidelberg, 2005. Springer Berlin Heidelberg.

\bibitem[AS15]{AbbeSandon2015a}
Emmanuel Abbe and Colin Sandon.
\newblock {Community Detection in General Stochastic Block Models: Fundamental
  Limits and Efficient Algorithms for Recovery}.
\newblock In {\em {Proceedings of the 2015 IEEE 56th Annual Symposium on
  Foundations of Computer Science (FOCS)}}, {FOCS '15}, pages 670--688,
  Washington, DC, USA, 2015. IEEE Computer Society.

\bibitem[ASW13]{azizyan2013minimax}
Martin Azizyan, Aarti Singh, and Larry Wasserman.
\newblock {Minimax theory for high-dimensional gaussian mixtures with sparse
  mean separation}.
\newblock In {\em {Advances in Neural Information Processing Systems}}, pages
  2139--2147, 2013.

\bibitem[AV07]{Kmeans++}
David Arthur and Sergei Vassilvitskii.
\newblock {K-means++: The Advantages of Careful Seeding}.
\newblock In {\em {Proceedings of the Eighteenth Annual ACM-SIAM Symposium on
  Discrete Algorithms}}, {SODA '07}, pages 1027--1035, Philadelphia, PA, USA,
  2007. Society for Industrial and Applied Mathematics.

\bibitem[BGL{\etalchar{+}}15]{cord}
F.~{Bunea}, C.~{Giraud}, X.~{Luo}, M.~{Royer}, and N.~{Verzelen}.
\newblock {Model Assisted Variable Clustering: Minimax-optimal Recovery and
  Algorithms}.
\newblock {\em ArXiv e-prints}, August 2015.

\bibitem[BGRV16]{pecok}
Florentina Bunea, Christophe Giraud, Martin Royer, and Nicolas Verzelen.
\newblock {PECOK: a convex optimization approach to variable clustering}.
\newblock {\em arXiv preprint arXiv:1606.05100}, 2016.

\bibitem[BRS18]{IsingBM}
Q.~{Berthet}, P.~{Rigollet}, and P.~{Srivastava}.
\newblock {Exact recovery in the Ising blockmodel}.
\newblock {\em Annals of Statistics (to appear)}, page arXiv:1612.03880, 2018.

\bibitem[CDF16]{2016arXiv160609190C}
S.~{Chr{\'e}tien}, C.~{Dombry}, and A.~{Faivre}.
\newblock {A Semi-Definite Programming approach to low dimensional embedding
  for unsupervised clustering}.
\newblock {\em ArXiv e-prints}, June 2016.

\bibitem[CGTS02]{kmedian}
Moses Charikar, Sudipto Guha, {\'E}va Tardos, and David~B. Shmoys.
\newblock {A Constant-Factor Approximation Algorithm for the k-Median Problem}.
\newblock {\em Journal of Computer and System Sciences}, 65(1):129--149, 2002.

\bibitem[CLX15]{2015arXiv151208425C}
Y.~{Chen}, X.~{Li}, and J.~{Xu}.
\newblock {Convexified Modularity Maximization for Degree-corrected Stochastic
  Block Models}.
\newblock {\em ArXiv e-prints}, December 2015.

\bibitem[CRV15]{Chin15}
Peter Chin, Anup Rao, and Van Vu.
\newblock {Stochastic Block Model and Community Detection in Sparse Graphs: A
  spectral algorithm with optimal rate of recovery}.
\newblock In Peter Gr{\"u}nwald, Elad Hazan, and Satyen Kale, editors, {\em
  {Proceedings of The 28th Conference on Learning Theory}}, volume~40 of {\em
  {Proceedings of Machine Learning Research}}, pages 391--423, Paris, France,
  03--06 Jul 2015. PMLR.

\bibitem[CX16]{chen2014statistical}
Yudong Chen and Jiaming Xu.
\newblock {Statistical-computational tradeoffs in planted problems and
  submatrix localization with a growing number of clusters and submatrices}.
\newblock {\em Journal of Machine Learning Research}, 17(27):1--57, 2016.

\bibitem[CY18]{ChenYang2018}
Xiaohui {Chen} and Yun {Yang}.
\newblock {Hanson-Wright inequality in Hilbert spaces with application to
  \$K\$-means clustering for non-Euclidean data}.
\newblock {\em arXiv e-prints}, page arXiv:1810.11180, Oct 2018.

\bibitem[DAM16]{DAM15}
Y.~Deshpande, E.~Abbe, and A.~Montanari.
\newblock {Asymptotic mutual information for the binary stochastic block
  model}.
\newblock In {\em {2016 IEEE International Symposium on Information Theory
  (ISIT)}}, pages 185--189, July 2016.

\bibitem[DKS17]{DBLP:journals/corr/abs-1711-07211}
Ilias Diakonikolas, Daniel~M. Kane, and Alistair Stewart.
\newblock {List-Decodable Robust Mean Estimation and Learning Mixtures of
  Spherical Gaussians}.
\newblock {\em CoRR}, abs/1711.07211, 2017.

\bibitem[FC17]{FeiChen2017}
Y.~{Fei} and Y.~{Chen}.
\newblock {Exponential error rates of SDP for block models: Beyond
  Grothendieck's inequality}.
\newblock {\em ArXiv e-prints}, 2017.

\bibitem[FC18]{FeiChen2018}
Y.~{Fei} and Y.~{Chen}.
\newblock {Hidden Integrality of SDP Relaxation for Sub-Gaussian Mixture
  Models}.
\newblock {\em ArXiv e-prints}, March 2018.

\bibitem[Gir15]{HDS}
Christophe Giraud.
\newblock {\em {Introduction to high-dimensional statistics}}, volume 139 of
  {\em {Monographs on Statistics and Applied Probability}}.
\newblock CRC Press, Boca Raton, FL, 2015.

\bibitem[GMZZ17]{Gao2017}
Chao Gao, Zongming Ma, Anderson~Y. Zhang, and Harrison~H. Zhou.
\newblock {Achieving Optimal Misclassification Proportion in Stochastic Block
  Models}.
\newblock {\em J. Mach. Learn. Res.}, 18(1):1980--2024, January 2017.

\bibitem[GV14]{guedon2014community}
Olivier Gu{\'e}don and Roman Vershynin.
\newblock {Community detection in sparse networks via Grothendieck's
  inequality}.
\newblock {\em arXiv preprint arXiv:1411.4686}, 2014.

\bibitem[HL17]{DBLP:journals/corr/abs-1711-07454}
Samuel~B. Hopkins and Jerry Li.
\newblock {Mixture Models, Robustness, and Sum of Squares Proofs}.
\newblock {\em CoRR}, abs/1711.07454, 2017.

\bibitem[HLL83]{holland1983stochastic}
Paul~W Holland, Kathryn~Blackmond Laskey, and Samuel Leinhardt.
\newblock {Stochastic blockmodels: First steps}.
\newblock {\em Social networks}, 5(2):109--137, 1983.

\bibitem[Hoe63]{Hoeffding}
Wassily Hoeffding.
\newblock {Probability Inequalities for Sums of Bounded Random Variables}.
\newblock {\em Journal of the American Statistical Association},
  58(301):13--30, 1963.

\bibitem[HWX16]{2016arXiv160206410H}
B.~{Hajek}, Y.~{Wu}, and J.~{Xu}.
\newblock {Semidefinite Programs for Exact Recovery of a Hidden Community}.
\newblock {\em ArXiv e-prints}, February 2016.

\bibitem[IMPV15]{iguchi2015tightness}
Takayuki Iguchi, Dustin~G Mixon, Jesse Peterson, and Soledad Villar.
\newblock {On the tightness of an SDP relaxation of k-means}.
\newblock {\em arXiv preprint arXiv:1505.04778}, 2015.

\bibitem[JMRT16]{JavanmardE2218}
Adel Javanmard, Andrea Montanari, and Federico Ricci-Tersenghi.
\newblock {Phase transitions in semidefinite relaxations}.
\newblock {\em Proceedings of the National Academy of Sciences},
  113(16):E2218--E2223, 2016.

\bibitem[KS17]{DBLP:journals/corr/abs-1711-07465}
Pravesh~K. Kothari and Jacob Steinhardt.
\newblock {Better Agnostic Clustering Via Relaxed Tensor Norms}.
\newblock {\em CoRR}, abs/1711.07465, 2017.

\bibitem[LCX18]{li2018convex}
Xiaodong Li, Yudong Chen, and Jiaming Xu.
\newblock {Convex Relaxation Methods for Community Detection}.
\newblock {\em arXiv preprint arXiv:1810.00315}, 2018.

\bibitem[LLL{\etalchar{+}}17]{li2017birds}
Xiaodong Li, Yang Li, Shuyang Ling, Thomas Strohmer, and Ke~Wei.
\newblock {When do birds of a feather flock together? k-means, proximity, and
  conic programming}.
\newblock {\em Mathematical Programming}, pages 1--47, 2017.

\bibitem[Llo82]{Lloyd}
S.~Lloyd.
\newblock {Least Squares Quantization in PCM}.
\newblock {\em IEEE Trans. Inf. Theor.}, 28(2):129--137, September 1982.

\bibitem[LR15]{LeiRinaldo}
Jing Lei and Alessandro Rinaldo.
\newblock {Consistency of spectral clustering in stochastic block models}.
\newblock {\em Ann. Statist.}, 43(1):215--237, 2015.

\bibitem[LZ16]{LuZhou2016}
Y.~{Lu} and H.~H. {Zhou}.
\newblock {Statistical and Computational Guarantees of Lloyd's Algorithm and
  its Variants}.
\newblock {\em ArXiv e-prints}, December 2016.

\bibitem[Moo17]{Moore2017}
Cristopher Moore.
\newblock {The Computer Science and Physics of Community Detection: Landscapes,
  Phase Transitions, and Hardness}.
\newblock {\em CoRR}, abs/1702.00467, 2017.

\bibitem[MPW16]{moitra2016robust}
Ankur Moitra, William Perry, and Alexander~S Wein.
\newblock {How robust are reconstruction thresholds for community detection?}
\newblock In {\em {Proceedings of the forty-eighth annual ACM symposium on
  Theory of Computing}}, pages 828--841. ACM, 2016.

\bibitem[MVW17]{mixon2016}
Dustin~G Mixon, Soledad Villar, and Rachel Ward.
\newblock {Clustering subgaussian mixtures by semidefinite programming}.
\newblock {\em Information and Inference: A Journal of the IMA}, 6(4):389--415,
  2017.

\bibitem[{Nda}18]{Ndaoud2018}
Mohamed {Ndaoud}.
\newblock {Sharp optimal recovery in the Two Component Gaussian Mixture Model}.
\newblock {\em arXiv e-prints}, page arXiv:1812.08078, Dec 2018.

\bibitem[PW07]{PengWei07}
Jiming Peng and Yu~Wei.
\newblock {Approximating K-means-type Clustering via Semidefinite Programming}.
\newblock {\em SIAM J. on Optimization}, 18(1):186--205, February 2007.

\bibitem[PW15]{2015arXiv150705605P}
A.~{Perry} and A.~S. {Wein}.
\newblock {A semidefinite program for unbalanced multisection in the stochastic
  block model}.
\newblock {\em ArXiv e-prints}, July 2015.

\bibitem[Roy17]{MartinNIPS}
M.~Royer.
\newblock {Adaptive Clustering through Semidefinite Programming}.
\newblock {\em Advances in Neural Information Processing Systems (NIPS)}, 2017.

\bibitem[RV13]{rudelson2013hanson}
Mark Rudelson and Roman Vershynin.
\newblock {Hanson-Wright inequality and sub-gaussian concentration}.
\newblock {\em Electron. Commun. Probab}, 18(82):1--9, 2013.

\bibitem[RV17]{Regev2017}
O.~{Regev} and A.~{Vijayaraghavan}.
\newblock {On Learning Mixtures of Well-Separated Gaussians}.
\newblock In {\em {2017 IEEE 58th Annual Symposium on Foundations of Computer
  Science (FOCS)}}, pages 85--96, Oct 2017.

\bibitem[VW04]{VEMPALA2004}
Santosh Vempala and Grant Wang.
\newblock {A spectral algorithm for learning mixture models}.
\newblock {\em Journal of Computer and System Sciences}, 68(4):841--860, 2004.
\newblock Special Issue on FOCS 2002.

\bibitem[YP14]{YunP14b}
Se{-}Young Yun and Alexandre Prouti{\`e}re.
\newblock {Accurate Community Detection in the Stochastic Block Model via
  Spectral Algorithms}.
\newblock {\em CoRR}, abs/1412.7335, 2014.

\end{thebibliography}

\appendix

\section{About CH-divergence}\label{section:CH}
\begin{lem}\label{lem:CH}
The CH-divergence defined by (\ref{eq:CH}) fulfills
$$D_{+}(q||p)  \leq {1\over 4\rho} \sum_{x} {(p_{x}-q_{x})^2\over p_{x}},\quad \textrm{when}\ \ \min_{x} {q_{x}\over p_{x}} \geq \rho >0.$$
\end{lem}
\noindent{\bf Proof.}
The function $f_{t}(y)=1-t+ty-y^t$ is convex and fulfills $f_{t}(1)=0=f'_{t}(1)$ and
$$f''_{t}(1+u)\leq {u^2\over 4 \rho},\quad \textrm{for all}\ \ 1+u\geq \rho\ \ \textrm{and}\ \ t\in[0,1].$$
Setting $u_{x}=(q_{x}-p_{x})/p_{x}$, we get the claimed result.

\section{Lower bound on misclassification probability in supervised Gaussian classification with unknown means}\label{sec:bayes}
In this section, we derive a lower-bound on the misclassification probability of the Bayes classifier, in the  Gaussian supervised classification problem with  two balanced classes, with identical spherical covariances $\Sigma_{k}=\sigma I_p$ and opposite means $\mu_{-1}= -\mu_1$ uniformly distributed on the Euclidean sphere $\partial B(0,\Delta/2)$ in $\R^p$.

We denote by $\cL=(X_{a},Z_{a})_{a=1,\ldots,n}$ the learning sample distributed as follows. The labels $Z_{1},\ldots,Z_{n}$ are i.i.d.\ with uniform distribution on $\ac{-1,1}$, a random vector $\mu\in\R^p$ is sampled  uniformly over the sphere $\partial B(0,\Delta/2)$ independently of  $Z_{1},\ldots,Z_{n}$, and, conditionally on $Z_{1},\ldots,Z_{n},\mu$, the $X_{a}$ are independent Gaussian random variables with mean $Z_{a}\mu$ and covariance $\sigma^2I_{p}$.

The classifier minimizing the misclassification probability $\P\cro{Z_{\rm new}\neq \widehat h(X_{\rm new})}$ over all the $\sigma(\cL)$-measurable  classifiers $\widehat h$ is the Bayes classifier  given by
$$\widehat h(x)=\textrm{sign}\big(\P\cro{Z=1|X=x,\cL}-\P\cro{Z=-1|X=x,\cL}\big).$$
Let us compute the Bayes classifier in our setting.  Indeed, the classification problem is scaling-invariant. 

For $\delta \in\ac{-1,1}$ and $x\in\R^p$, we have
$$\P\cro{Z=\delta |X=x,\cL}=\int_{\partial B(0,\Delta/2)} \P[Z=\delta | X=x,\cL,\mu]\,d\P[\mu | X=x, \cL].$$
Direct computations give
$$\P[Z=\delta | X=x,\cL,\mu]=\P[Z=\delta | X=x,\mu]={e^{-0.5\|\delta x-\mu\|^2/\sigma^2}\over e^{-0.5\|x+\mu\|^2/\sigma^2}+e^{-0.5\|x-\mu\|^2/\sigma^2}}$$
and
$$d\P[\mu | X=x, \cL] \propto \pa{e^{-0.5\|x+\mu\|^2/\sigma^2}+e^{-0.5\|x-\mu\|^2/\sigma^2}} e^{-0.5 \sum_{a}\|Z_{a}X_{a}-\mu\|^2/\sigma^2}.$$
Hence, by using that $\|\mu\|=\Delta/2$ on $\partial B(0,\Delta/2)$, and by denoting by $\gamma$ the uniform distribution on $\partial B(0,\Delta/2)$,  we obtain
\begin{align*}
\P\cro{Z=\delta |X=x,\cL}&={\int_{\partial B(0,\Delta/2)}e^{-0.5\|\delta x-\mu\|^2/\sigma^2}e^{-0.5 \sum_{a}\|Z_{a}X_{a}-\mu\|^2/\sigma^2}d\gamma(\mu)\over \int_{\partial B(0,\Delta/2)}\pa{e^{-0.5\|x+\mu'\|^2/\sigma^2}+e^{-0.5\|x-\mu'\|^2/\sigma^2}}e^{-0.5 \sum_{a}\|Z_{a}X_{a}-\mu'\|^2/\sigma^2}d\gamma(\mu')}\\
&={\int_{\partial B(0,\Delta/2)}   e^{-\langle \delta x+\sum_{a}Z_{a}X_{a},\mu\rangle/\sigma^2}d\gamma(\mu)
\over \int_{\partial B(0,\Delta/2)}e^{-\langle x+\sum_{a}Z_{a}X_{a},\mu'\rangle/\sigma^2}d\gamma(\mu')+ \int_{\partial B(0,\Delta/2)}e^{-\langle -x+\sum_{a}Z_{a}X_{a},\mu'\rangle/\sigma^2}d\gamma(\mu')}\,.
\end{align*}
Since $F(v)=\int_{\partial B(0,\Delta/2)} e^{\langle v,\mu\rangle} d\gamma(\mu)$ depends only on $\|v\|$ and is monotone increasing with $\|v\|$, we obtain that 
\begin{align*}
\P\cro{Z=1|X=x,\cL}>\P\cro{Z=-1|X=x,\cL} & \ \Longleftrightarrow \ \big\|x+\sum_{a}Z_{a}X_{a}\big\|^2 >\big\|-x+\sum_{a}Z_{a}X_{a}\big\|^2\\
&\ \Longleftrightarrow \  \big\langle x,\sum_{a}Z_{a}X_{a}\big\rangle >0,
\end{align*}
and finally
$$\widehat{h}(x)= \textrm{sign}\pa{\bigg\langle {1\over n}\sum_{a=1}^n Z_{a}X_{a},x\bigg\rangle}.$$

For any $\sigma>0$, the probability of misclassification of the Bayes classifier is given by \begin{align*}
\P\cro{Z_{\rm new}\neq \widehat h(X_{\rm new})} &=  \int_{\partial B(0,\Delta/2)} \P\cro{Z\widehat h(X)<0\Big|\mu}\, d\gamma(\mu)\\
&=  \int_{\partial B(0,\Delta/2)} \P\cro{\big\langle \mu + {\sigma\over \sqrt{n}} \epsilon, \mu+\sigma\epsilon'\big\rangle<0\Big|\mu}\, d\gamma(\mu),
\end{align*}
where $\epsilon$ and $\epsilon'$ are two independent standard Gaussian random variables in $\R^p$.
The above conditional probability is invariant over $\partial B(0,\Delta/2)$ hence we only need to evaluate it for a fixed $\mu \in \partial B(0,\Delta/2)$, say $\mu_{\Delta}=[\Delta/2,0,\ldots,0]$. Let us set $W=-2\pa{\Delta \sqrt{1+1/n}}^{-1}\langle \mu_{\Delta},  {1\over \sqrt{n}}\epsilon+\epsilon'\rangle$ which follows a standard Gaussian distribution in $\R$ and  $S=-\langle \epsilon,\epsilon'\rangle$. Then, we have
\begin{align*}
\P\cro{Z_{\rm new}\neq \widehat h(X_{\rm new})} &=\P\cro{\big\langle \mu_{\Delta} + {\sigma\over \sqrt{n}} \epsilon, \mu_{\Delta}+\sigma\epsilon'\big\rangle<0}\\
&= \P\cro{{\Delta^2\over 4\sigma^2} < {\Delta\over 2\sigma}\sqrt{1+{1\over n}}\,W+{1\over \sqrt{n}}S}\ .
\end{align*}
We observe that $(-W,S)$ has the same distribution as $(W,S)$, hence by a union bound
$$ \P\cro{{\Delta^2\over 4\sigma^2} < {\Delta\over 2\sigma}\sqrt{1+{1\over n}}\,|W|+{1\over \sqrt{n}}S} \leq 2  \P\cro{{\Delta^2\over 4\sigma^2} < {\Delta\over 2\sigma}\sqrt{1+{1\over n}}\,W+{1\over \sqrt{n}}S}.$$
By using that the distributions of $S$ and $W$ are symmetric and that $a\vee (2b-2a) \geq a \vee (b/2)$, we get
\begin{align*}
 \P&\cro{{\Delta^2\over 4\sigma^2} < {\Delta\over 2\sigma}\sqrt{1+{1\over n}}\,|W|+{1\over \sqrt{n}}S}\\ 
 &\geq  \P\cro{{\Delta^2\over 4\sigma^2} < {1\over \sqrt{n}}S}+ \P\cro{{\Delta^2\over 2\sigma^2} < {\Delta\over 2\sigma}\sqrt{1+{1\over n}}\,|W|\,; {1\over \sqrt{n}}|S|\leq {\Delta^2\over 4\sigma^2}}\\
& \geq  \P\cro{{\Delta^2\over 4\sigma^2} < {1\over \sqrt{n}}S} \bigvee \pa{\P\cro{{\Delta\over \sigma} < |W|}-\P\cro{{1\over \sqrt{n}}|S|> {\Delta^2\over 4\sigma^2}}}\\
& \geq  \P\cro{{\Delta^2\over 4\sigma^2} < {1\over \sqrt{n}}S} \bigvee \pa{2\P\cro{{\Delta\over \sigma} < W}-2\P\cro{{1\over \sqrt{n}}S> {\Delta^2\over 4\sigma^2}}}\\
&\geq \P\cro{{\Delta^2\over 4\sigma^2} < {1\over \sqrt{n}}S} \bigvee {1\over 2}\P\cro{{\Delta\over \sigma} < W}.
\end{align*}
Putting pieces together, we get 
$$4\P\cro{Z_{\rm new}\neq \widehat h(X_{\rm new})} \geq  \P\cro{{\Delta\over \sigma} < W} \bigvee \P\cro{{\Delta^2\over 4\sigma^2} < {1\over \sqrt{n}}S}.$$

We observe that $S$ is distributed as the product of a standard Gaussian real random variable with the square-root of an independent $\chi^2$ random variable with $p$ degrees of freedom. Since a $\chi^2$ random variable with $p$ degrees of freedom is larger than $p/2$ with probability larger than $1/2$, we then have
$$\P\cro{S>{\sqrt{n}\,\Delta^2\over 4\sigma^2}}\geq {1\over 2} \P\cro{\sqrt{p\over 2}\, W >{\sqrt{n}\Delta^2\over 4\sigma^2}} \geq {1\over 2} \P\cro{W > {\sqrt{n}\,\Delta^2\over \sqrt{8p}\,\sigma^2}}.$$
We then obtain the lower bound on the Bayes probability of misclassification
$$\P\cro{Z_{\rm new}\neq \widehat h(X_{\rm new})}\geq c\exp\pa{-c'\pa{{\Delta^2\over \sigma^2}\wedge {n \Delta^4 \over p \sigma^4}}},$$
for some numerical constants $c,c'>0$.

\end{document}